\documentclass[11pt]{article}

\usepackage{amssymb,commath,
graphics,amsmath,amsthm,amsopn,amstext,amsfonts,hyperref,fullpage,graphicx,
fancybox,verbatim,cite}

\usepackage{fancyhdr,color,mathtools, algorithmic,cite}
\usepackage{xcolor}
\usepackage{verbatim}

\DeclareFontFamily{U}{mathx}{\hyphenchar\font45}
\DeclareFontShape{U}{mathx}{m}{n}{
      <5> <6> <7> <8> <9> <10>
      <10.95> <12> <14.4> <17.28> <20.74> <24.88>
      mathx10
      }{}
\DeclareSymbolFont{mathx}{U}{mathx}{m}{n}
\DeclareFontSubstitution{U}{mathx}{m}{n}
\DeclareMathAccent{\widecheck}{0}{mathx}{"71}

\newtheorem{theorem}{Theorem}
\newtheorem{proposition}[theorem]{Proposition}
\newtheorem{lemma}[theorem]{Lemma}
\newtheorem{corollary}[theorem]{Corollary}
\newtheorem{definition}[theorem]{Definition}
\newtheorem{example}[theorem]{Example}

\newcommand{\gap}{\vspace{0.1in}}

\newcommand{\epc}{\hspace{1pc}}
\newcommand{\thalf}{{\textstyle{\frac{1}{2}}}}

\newcommand{\onebld}{{\bf 1}}
\newcommand{\wt}{\widetilde}
\newcommand{\wh}{\widehat}
\newcommand{\ball}{\mathbb{B}}
\newcommand{\dist}{\mbox{dist}}

\usepackage{tikz}
\usetikzlibrary{patterns, arrows.meta}

\title{The Minimization of Piecewise Functions: Pseudo Stationarity \\ [0.2in]
{\small \sl This paper is dedicated to Professor Roger J.B.\ Wets on the occasion
of his 85th birthday, \\ [-0.1in]
and for his pioneering research on the subject of our work.}}

\author{
Ying Cui\footnote{Department of Industrial and Systems Engineering, University
of Minnesota, Minneapolis, U.S.A.\ 55455. {\tt Email: yingcui@umn.edu.}} \and
Junyi Liu\footnote{Department of Industrial Engineering, Tsinghua University,
Beijing, China 100084.  {\tt Email: junyiliu@mail.tsinghua.edu.cn.}}
\and Jong-Shi Pang\footnote{The
Daniel J.\ Epstein Department of Industrial and Systems Engineering, University of
Southern California, Los Angeles, U.S.A.\ 90089.
This work was based on research supported by the  U.S.\ Air Force Office of
Sponsored Research under grants FA9550-18-1-0382 and FA9550-22-1-0045. {\tt Email: jongship@usc.edu.}}
}

\date{Original: January 30, 2022; Revised: July 1, 2022}

\begin{document}

\maketitle

\begin{abstract}
\noindent There are many significant applied contexts that require the solution of discontinuous
optimization problems in finite dimensions.  Yet these problems are very difficult, both
computationally and analytically.  With the functions being discontinuous and a minimizer (local or global)
of the problems, even if it exists, being impossible to verifiably compute, a foremost
question is what kind of ``stationary solutions'' one can expect to obtain; these solutions
provide promising candidates for minimizers; i.e., their defining conditions are necessary for optimality.
Motivated by recent results on sparse optimization, we introduce in this paper such a kind of solution,
termed ``pseudo B- (for Bouligand) stationary solution'',
for a broad class of discontinuous
optimization problems with objective and constraint defined by indicator functions of the positive real
axis composite with functions that are possibly nonsmooth.  We present two approaches for computing such a solution.
One approach is based on lifting the problem to a higher
dimension via the epigraphical formulation of the indicator functions; this requires the addition of some
auxiliary variables.  The other approach is based on certain continuous (albeit not necessarily differentiable)
piecewise approximations of the
indicator functions and the convergence to a pseudo B-stationary solution of the
original problem is established.  The conditions for convergence are discussed and illustrated by an example.
\end{abstract}

\section{Introduction}

There are many significant applied contexts that require the solution of discontinuous
optimization problems in finite dimensions.  Yet these problems are very difficult, both
computationally and analytically.  Like all nonconvex problems, the computational task of a minimizer, local
or global, is prohibitively challenging, if not impossible.  Thus a realistic goal is to try to compute a solution
that is a promising candidate for a minimizer, in particular, one that satisfies
some necessary conditions of a local minimizer, i.e., a stationary solution of some sort.  There is
a very long tradition of investigation of such conditions, starting from the classical one of a zero objective gradient
for an unconstrained differentiable optimization problem to the use of some
advanced subdifferentials from variational analysis \cite{RockafellarWets09} to obtain a set-inclusion problem.
With the functions involved being discontinuous, an insightful understanding, constructive characterization,
and profitable employment of the latter subdifferentials all are not easy tasks, especially when there are
constraints that are embedded in the objective function taking on infinite value.

\gap

As an alternative, the
idea of approximating the discontinuous functions by smooth, say continuously differentiable, functions appears
promising.
 {One of the first smoothing methods for optimization of discontinuous functions was
proposed in \cite{GupalNorkin77}.  In general, such a method involves smoothing a locally integrable function
by mollifier approximations using integration.}
A sequence of the resulting smoothed problems is then solved and their stationary solutions are used to define
a stationarity concept for the original problem.  Studied comprehensively in the seminal paper \cite{ENWets95}
under the framework of minimizing an extended-valued function,
this approach leads to the definition of ``mollifier subgradients''  that serve as the target of a
computational resolution for a discontinuous optimization problem.  See \cite{JongenStein-I,JongenStein-II} for the
applications of mollifier induced smoothing functions, called ``averaged function'' in \cite[Definition~3.1]{ENWets95}, to
nonlinear and semi-infinite programming, respectively.
The paper \cite{Chen12} 
gives a survey of smoothing methods for nonconvex, nonsmooth optimization problems
with the underlying functions being continuous.
It is important to point out while the theory in \cite{ENWets95,Chen12} is quite general, the applied problems surveyed
in \cite{Chen12} are all of the kind of a univariate nonsmooth function composite with a smooth function; for these
composite functions, smoothing
is applied to the univariate component while the composition is maintained.  Thus it is fair to say that to date, the practical aspect of
smoothing is restricted to this class of composite functions where the convolutional operation requires only the evaluation of
integrals of scalar functions of one real variable; the approach becomes most effective when the resulting averaged functions
are explicitly available for
computational purposes with the convolution operation staying in the background for general analysis. {A related paper \cite{EN98} discusses smoothing by introducing artificial random variables that
also allow the use of stochastic gradient methods for solving the approximated problems.}

\gap

Supported by diverse source problems, this paper studies a broad class of discontinuous optimization problems with objective
and constraints defined by piecewise functions modeled as the products of nonconvex
nonsmooth functions and the indicator functions of the positive/nonnegative real axis, called {\sl Heaviside functions}. {Invented by Oliver Heaviside (1850--1925) in his pioneering work on differential equations for the
study of electromagnetic waves \cite{Heaviside88}, the ``open'' Heaviside function is the indicator of the open interval
$( \, 0,\infty \, )$.  For our purpose, we include the indicator of the closed interval
$[ \, 0,\infty \, )$ also as a Heaviside function.}
With the advance of nonconvex nonsmooth optimization as documented in \cite{CuiPang2021}, our goal in approximating these problems is to least
disrupt the nonsmoothness and nonconvexity of the given functions and do so only when needed, for instance in designing practical
computational methods. {Thus, unlike the previous work \cite{GupalNorkin77,Batukhtin93,ENWets95,EN98,EN04} that smooth the discontinuous Heaviside function,
the approximating problems developed in the present paper remain nonsmooth and nonconvex but are solvable by various surrogation methods;
see \cite[Chapter~7]{CuiPang2021}.}
In addition to this departure from the smooth approximations, we adopt the approach to 
define an approximation-independent and subdifferential-free stationarity
condition, called {\sl pseudo B- (for Bouligand) stationarity}, that is necessarily satisfied by a local minimizer of the given
discontinuous optimization problem.
The definition is motivated by a recent study \cite{GomezHePang21} of the $\ell_0$-optimization problem that lies at the center of sparse
optimization in statistical estimation \cite{HastieTibWeinwright15}.  The well-known $\ell_0$-function is defined by:
$| \, t \, |_0 \triangleq \left\{ \begin{array}{ll}
1 & \mbox{if $t \neq 0$} \\
0 & \mbox{otherwise.}
\end{array} \right.$
We describe two constructive approaches for computing/approximating such a pseudo B-stationary solution.
Omitting the details that are left for a follow-up
algorithmic study, we emphasize that these approaches can be computationally implemented in practice
by difference-of-convex programming based algorithms \cite{LeThiPham05,PhamLeThi97,LuZhouSun19,PangRazaAlvarado16,QiCuiLiuPang19}
when the involved functions are of this kind, and more generally, by surrogation methods
\cite[Chapter~7]{CuiPang2021} for broad classes of nonsmooth nonconvex optimization
problems.

\gap

The rest of the paper is organized as follows.  The next section begins with the formal definition of the problem to be studied followed by
some preliminary remarks.   Section~\ref{sec:source problems} presents a host of discontinuous piecewise
functions arising from various optimization contexts that are unified by our central problem.  A main result in
Section~\ref{sec:global minimizers} identifies a principal sign condition that plays a central role throughout the paper;
this result connects our problem with two related problems in the literature in terms of their global minima.  As a remedy to
the computational intractability
of these global minima, Section~\ref{sec:pseudo} defines a pseudo B-stationarity concept that is amenable to computation.  Two
constructive approaches for computing such a stationary solution is presented in Sections~\ref{sec:epi} and \ref{sec:approximation computation}, respectively.
In particular, the approximation approach described in the last section can be traced back to a pioneering paper by Roger Wets and his collaborators that we
expand in Section~\ref{sec:approximations Heaviside} and to whom we dedicate our work.

\section{Problem Definition and Preliminary Discussion}

Consider the following piecewise 
optimization problem:
\begin{equation} \label{eq:original POP}
\begin{array}{ll}
\displaystyle{
\operatornamewithlimits{\mbox{\bf minimize}}_{x \in X}
} & \Phi(x) \, \triangleq \, c(x) +
\displaystyle{
\sum_{k=1}^K
} \, \varphi_k(x) \, \onebld_{( \, 0,\infty \, )}( g_k(x)) \\ [0.15in]
\mbox{\bf subject to} & \displaystyle{
\sum_{\ell=1}^L
} \, \phi_{\ell}(x) \, \onebld_{( \, 0,\infty \, )}( h_{\ell}(x)) \, \leq \, b, \epc \mbox{called the {\sl functional constraint}},
\end{array}
\end{equation}
where $\onebld_{( 0,\infty )}$ is the ``open'' Heaviside function given by
\[
\onebld_{( \, 0, \infty \, )}(s) \, \triangleq \, \left\{ \begin{array}{ll}
1 & \mbox{if $s \in ( \, 0, \infty \, )$} \\ [5pt]
0 & \mbox{otherwise}
\end{array} \right.
\]
and the following holds [the polyhedrality of $X$ is not needed in several results; it is stated as a blanket
assumption primarily to avoid the use of advanced constraint qualifications in the context of tangent cones]:

\gap

\fbox{
\parbox{6.5in}{
{\bf Blanket assumption:} $X$ is a polyhedron  {contained in the open subset ${\cal O}$ of $\mathbb{R}^n$} and
$c$, $\{ \, \varphi_k,g_k \, \}_{k=1}^K$, and $\{ \, \phi_{\ell},h_{\ell} \}_{\ell=1}^L$
are B-differentiable (where B is for Bouligand) functions from ${\cal O} \to \mathbb{R}$. \hfill $\Box$
}}

\gap

By definition, a function $\psi : {\cal O} \to \mathbb{R}$ is B-differentiable \cite[Definition~4.1.1]{CuiPang2021}
at $\bar{x} \in {\cal O}$ if $\psi$ is locally Lipschitz continuous at $\bar{x}$ (i.e., Lipschitz continuous in an open neighborhood of $\bar{x}$)
and directionally differentiable there; i.e., the elementary one-sided directional derivative
\[
\psi^{\, \prime}(\bar{x};v) \, \triangleq \, \displaystyle{
\lim_{\tau \downarrow 0}
} \, \displaystyle{
\frac{\psi(\bar{x} + \tau v) - \psi(\bar{x})}{\tau}
}
\]
exists for all $v \in \mathbb{R}^n$.  We make several immediate remarks about the formulation (\ref{eq:original POP}); foremost
is to note that the functional constraint adds considerable challenges to this problem, without which the
analysis simplifies somewhat.  Other remarks are as follows:

\gap

$\bullet $  While the open Heaviside function ${\bf 1}_{(\,0, \,\infty)}(\bullet)$
is lower semicontinuous on the real line, a product such as $\varphi_k( \bullet ) \, \onebld_{( \, 0,\infty \, )}( g_k( \bullet ))$ may not be
lower semicontinuous at $\bar{x} \in g_k^{-1}(0)$ unless $\varphi_k$ is nonnegative there.  For the
problem (\ref{eq:original POP}), this sign condition on the pairs of functions $\{ \varphi_k,g_k \}_{k=1}^K$ and $\{ \phi_{\ell},h_{\ell} \}_{\ell=1}^L$
will persist throughout the paper.  This restriction is responsible for the closedness of the feasible region and for the existence of
 minimizers of the problem, and thus for their relaxations, such as that of a pseudo B-stationary point to be defined later.

\gap

$\bullet $ The well-known $\ell_0$-function $| \, \bullet \, |_0$ in sparsity estimation \cite{HastieTibWeinwright15}
can be written as:
\begin{equation} \label{eq:ell0 fnc}
| \, s \, |_0 \, = \, \onebld_{( \, 0, \infty \, )}(s) + \onebld_{( \, 0, \infty \, )}(-s) \, = \, \onebld_{( \, 0, \infty \, )}(| s | ).
\end{equation}
This is a prominent applied instance of the open Heaviside function and provides much motivation for the theory developed in this paper.
Conversely, the developed theory also provides a deeper understanding of the sparse optimization problem
 {as a result of (\ref{eq:ell0 func}), any approximation of the open Heaviside
function can immediately be specialized to the $\ell_0$-function.}

\gap

$\bullet $ A term $\psi(x) \, \onebld_{[ \, 0, \infty \, )}(f(x))$ involving the ``closed'' Heaviside function
$\onebld_{[ \, 0, \infty \, )}(s) \, \triangleq \, \left\{ \begin{array}{ll}
1 & \mbox{if $s \in [ \, 0, \infty \, )$} \\
0 & \mbox{otherwise}
\end{array} \right.$ can be written as:
\[
\psi(x) \, \onebld_{[ \, 0, \infty \, )}(f(x)) \, = \, \psi(x) - \psi(x) \, \onebld_{( \, 0, \infty \, )}(-f(x)).
\]
So the formulation (\ref{eq:original POP}) encompasses products of this kind that involves the closed Heaviside function.
Subsequently, conditions imposed on the functions in (\ref{eq:original POP}) easily translate to corresponding conditions
for functions involving the closed Heaviside function via the above relation.

\gap

In general, a locally Lipschitz function does not need to be directionally differentiable; however, for a univariate function
$f : ( -\ell,u ) \to \mathbb{R}$ defined on an interval, if $f$ is locally Lipschitz and monotone (i.e., nondecreasing or nonincreasing),
then the one-sided derivatives:
\begin{equation} \label{eq:dd of univariate}
f^{\, \prime}(t;\pm 1) \, \triangleq \, \displaystyle{
\lim_{\tau \downarrow 0}
} \, \displaystyle{
\frac{f(t \pm \tau) - f(t)}{\tau}
}, \epc t \, \in \, ( -\ell,u )
\end{equation}
exist with signs determined by the monotonicity.  For an integer $N > 0$, we let $[ N ] \triangleq \{ 1, \cdots, N \}$.

\section{Some Source Problems} \label{sec:source problems}

The product functions in (\ref{eq:original POP}) are discontinuous piecewise functions; they include several interesting special
cases, which we highlight below.  These cases illustrate the versatility of the Heaviside functions in modeling
a host of discontinuous functions in diverse contexts.

\gap

$\bullet $ {\sl Cost-efficient variable selection:} This is an extension of the sparsity optimization problem in statistical
estimation \cite{HastieTibWeinwright15} in which there is a cost associated the collection of data in the modeling process
\cite{BianChen20,YuFuLiu21,Yue10}
which we may formulate as either a soft penalty embedded in the objective function:
\begin{equation} \label{eq:cost-sensitive variable selection}
\displaystyle{
\operatornamewithlimits{\mbox{\bf minimize}}_{x \in X \subseteq \mathbb{R}^n}
} \ e(x) + \lambda \, \displaystyle{
\sum_{i=1}^n
} \, c_i \, | \, x_i \, |_0,
\end{equation}
or as a budget-type constraint:
\begin{equation} \label{eq:constraint form cost-sensitive variable selection}
\displaystyle{
\operatornamewithlimits{\mbox{\bf minimize}}_{x \in X \subseteq \mathbb{R}^n}
} \ e(x)
\epc \mbox{\bf subject to } \ \displaystyle{
\sum_{i=1}^n
} \, c_i \, | \, x_i \, |_0 \, \leq \, b,
\end{equation}
where $e(x)$ is a loss function, $\lambda > 0$ is a given parameter, and the coefficients $c_i$ {and the right-hand constant $b$} are positive.   {In addition to the application in
modern-day sparse optimization, the term $c_i | x_i |_0$ is a
classical modeling device in operations research applications as a set-up cost of an activity.  Namely, a cost is
incurred when there is a nonzero level of
the $i$th activity and zero otherwise.  The functional constraint expresses the available budget for the set-up or variable-selection
costs applicable in both the data or operational contexts.
A cost associated with the level of the activity is included in either the objective term $e(x)$ or the constraint $X$.}
One mathematical feature of the summation term in both formulations is worth noting: namely, all the coefficients associated
with the $\ell_0$-functions are positive, thus, satisfying the sign condition mentioned before.  These
problems suggest that the multiplicative functions
$\varphi_k(x)$ and $\phi_{\ell}(x)$ in the general problem (\ref{eq:original POP}) may be interpreted as costs in some applied models.

\gap

$\bullet $ {\sl Piecewise functions on complementary regions:}  Consider a simple example of a piecewise function which may or
may not be continuous:
\begin{equation} \label{eq:deterministic 3-piece}
\Psi(x) \, = \, \left\{ \begin{array}{ll}
\psi_1(x) & \mbox{if $a \leq f(x) \leq b$} \\ [5pt]
\psi_2(x) & \mbox{if $f(x) < a$} \\ [5pt]
\psi_3(x) & \mbox{if $f(x) > b$},
\end{array} \right.
\end{equation}
for some scalars $a$ and $b$ satisfying $-\infty \leq a < b \leq \infty$.  We then have
\[ \begin{array}{rl}
\Psi(x) \, = & \psi_1(x) \, \onebld_{[ \, 0,\infty \, )}\left( \, \min( \, b - f(x), f(x) - a \, ) \, \right) \ + \\ [0.1in]
& \psi_2(x) \, \onebld_{( \, 0,\infty \, )}( a - f(x) ) + \psi_3(x) \, \onebld_{( \, 0,\infty \, )}( f(x) - b ) \\ [0.1in]
= & \psi_1(x) - \psi_1(x) \, \onebld_{( \, 0,\infty \, )}\left( \, \max( \, f(x) - b, a - f(x) \, ) \, \right) \ + \\ [0.1in]
& \psi_2(x) \, \onebld_{( \, 0,\infty \, )}( a - f(x) ) + \psi_3(x) \, \onebld_{( \, 0,\infty \, )}( f(x) - b ),
\end{array}
\]
whose validity is regardless of the continuity of $\Psi$ on the sets $f^{-1}(a)$ and $f^{-1}(b)$.  An interesting application
of this class of piecewise function concerns the constrained optimization with (soft) penalty.  Specifically, let's say that
we wish to minimize a function $\psi_1(x)$ when the constraint $f(x) \geq 0$ is satisfied and there is a penalty $\psi_2(x)$
when the constraint is not satisfied.  This problem can be formulated as minimizing the combined objective:
$\psi_1(x) \, \onebld_{[ \, 0,\infty \, )}(f(x)) + \psi_2(x) \, \onebld_{( \, 0,\infty \, )}(-f(x))$.  A special case of
(\ref{eq:deterministic 3-piece}) is when the three component functions $\psi_i(x)$ are constants.  More generally, a quantized
function is a discontinuous step function that can very easily be described by a generalization of (\ref{eq:deterministic 3-piece})
with arbitrary (finite) number of mutually disjoint regions (intervals in the case of a univariate variable)
within each of which the overall function is a constant.  The references \cite{LongYinXin21,Xin19,YinZhangQiXin19} have employed
such quantized minimization problems for the training of deep neural networks.

\gap

$\bullet $ {\sl Best constraint selection:}
In contrast to the best variable selection problems (\ref{eq:cost-sensitive variable selection})
 {or} (\ref{eq:constraint form cost-sensitive variable selection}), and yet
so far at best minimally studied
in the literature, the best constraint selection is a generalization of best variable selection and can be modeled using
the summation: $\displaystyle{
\sum_{k=1}^K
} \, \onebld_{[ \, 0,\infty )}(f_k(x))$, where the family $\{ f_k(x) \geq 0 \}_{k=1}^K$ consists of the constraints to be selected.
A related problem is when there is a cost $c_k(x)$ associated with the constraint $f_k(x) \geq 0$ being disrupted.  This problem can
be formulated as minimizing the weighted sum: $\displaystyle{
\sum_{k=1}^K
} \, c_k(x) \, \onebld_{( \, 0,\infty )}(-f_k(x))$ or imposing a constraint defined by the sum.
Our work offers a pathway for the potential applications of best constraint selection in instances like these.

\gap

$\bullet $ {\sl On-off constraints by indicator variables:}  In the literature such as \cite{BLTrWiese15,HBCornuejols12},
such a constraint is $f_k(x) \geq 0$ if $y_k = 1$ where $y_k \in \{ 0,1 \}$ is a binary variable that appears only in the constraint and for
the sole purpose of turning on the constraint $f_k(x) \geq 0$.  Equivalently, such an on-off constraint is equivalent to:
$f_k(x) \geq 0$ if $y_k > 0$ with $y_k$ restricted to be a continuous variable in the interval $[ 0,1 ]$.  In turn, the latter constraint
is equivalent to $y_k f_k(x) \geq 0$ with $y_k \in [ 0,1 ]$, and instead of the constraint, a term $\onebld_{[ \, 0,\infty )}(y_k f_k(x))$ can
be added to the objective as part of the overall function to be minimized.

\gap

$\bullet $ {\sl Binary classifications by the sign function:}  In binary classification,
it is customary to use the sign function to separate
two classes.  For instance, given a binary number $\sigma = \pm 1$ denoting classes A and B, respectively and with $f(x)$ as the feature-dependent
classification function, we classify the outcome from this function as A if $f(x) > 0$ and $B$ if $f(x) < 0$.  This classification scheme can by modeled
by the composite indicator function: $\onebld_{( \, 0,\infty \, )}(-\sigma f(x))$, which aims to count the number of misclassifications among the observational data.
See \cite{QiCuiLiuPang19} for a recent application of this
formulation in the context of individualized decision making under uncertainty for medical treatment.  Classification with margin is an extension that allows for
minor errors; it classifies the outcome as A if $f(x) \geq \varepsilon$ and $B$ if $f(x) \leq -\varepsilon$,
where $\varepsilon > 0$ is a small margin.  This can be modeled by $\onebld_{( \, 0,\infty \, )}(\varepsilon - \sigma f(x))$,
which again counts the number of misclassifications.

\gap

$\bullet $ {\sl Products of indicators:}  It is clear that $\left( \, \onebld_{[ \, 0,\infty \, )}(f(x)) \, \right) \,
\left( \, \onebld_{[ \, 0,\infty \, )}(g(x)) \, \right) = \onebld_{[ \, 0,\infty \, )}(\min( f(x),g(x) ))$; thus products of closed Heaviside functions
can be combined into a single Heaviside function of the same kind via the use of the pointwise minimum operator.  Less clear is the product of a closed and an
open Heaviside function, say $\left( \, \onebld_{[ \, 0,\infty \, )}(f(x)) \, \right) \,
\left( \, \onebld_{( \, 0,\infty \, )}(g(x)) \, \right)$.  Nevertheless, we have
\[ \begin{array}{lll}
\left( \, \onebld_{[ \, 0,\infty \, )}(f(x)) \, \right) \,
\left( \, \onebld_{( \, 0,\infty \, )}(g(x)) \, \right) & = & \left( \, \onebld_{( \, 0,\infty \, )}(g(x)) \, \right) \left[ \,
1 - \left( \, \onebld_{( \, 0,\infty \, )}(-f(x)) \, \right) \, \right] \\ [0.2in]
& = & \onebld_{( \, 0,\infty \, )}(g(x)) - \onebld_{( \, 0,\infty \, )}(\min( -f(x),g(x) )),
\end{array} \]
which is the difference of two open Heaviside functions.  An example where a product of an open and a closed Heaviside function may occur is a modification of
the piecewise function (\ref{eq:deterministic 3-piece}): say $\psi(x) = \psi_1(x)$ if $a \leq f(x) < b$.  We see that
\[ \begin{array}{lll}
\psi(x) & = & \psi_1(x) \, \left( \, \onebld_{[ \, 0,\infty \, )}(f(x) - a) \, \right) \, \left( \, \onebld_{( \, 0,\infty \, )}( b - f(x) ) \, \right) \\ [0.1in]
& = & \psi_1(x) \, \left[ \, \onebld_{( \, 0,\infty \, )}(b - f(x)) - \onebld_{( \, 0,\infty \, )}(a - f(x)) \, \right].
\end{array} \]
$\bullet $ {\sl Probabilistic functions and conditional expectations:}
The use of the Heaviside functions in equivalent formulations of probabilistic functions is well known; see \cite{CuiLiuPang21} for
a recent comprehensive study of a nonconvex nonsmooth approach for chance-constrained stochastic programs that is built on this
fundamental formulation.    The simplest example is the probabilistic function $\mathbb{P}_{\tilde{z}}(\Psi(x,\tilde{z}) \geq 0)$,
for some bivariate function $\Psi(x,z)$.  In the cited reference, we have considered a broad class of such functions
defined as the difference of two convex functions each being the pointwise maximum of finitely many convex functions.  An abstraction
of such a difference-of-convex function $\Psi(\bullet,z)$ is a piecewise function of which the following bivariate extension of
(\ref{eq:deterministic 3-piece}) is an example:
%
\[
\Psi(x,z) \, \triangleq \, \left\{ \begin{array}{ll}
\psi_1(x,z) & \mbox{if $a \leq f(x,z) \leq b$} \\ [5pt]
\psi_2(x,z) & \mbox{if $f(x,z) < a$} \\ [5pt]
\psi_3(x,z) & \mbox{if $f(x,z) > b$}.
\end{array} \right.
\]
Appendix B of the paper \cite{YifanCui21} contains many piecewise functions of the above kind arising from individualized decision-making under
partial identification.  Such a piecewise random functional gives rise to an expectation of piecewise functions such as
$\mathbb{E}_{\tilde{z}}\left[ \, \Psi(x,\tilde{z}) \, \right]$ with $\Psi$ given above.

\gap

Additionally, consider the conditional expectation:
\[
\mathbb{E}_{\tilde{z}}\left[ \, \phi(x,\tilde{z}) \mid f(x,\tilde{z}) \leq 0 \, \right]
\, \triangleq \, \displaystyle{
\frac{\mathbb{E}_{\tilde{z}}\left[ \, \phi(x,\tilde{z}) \, \onebld_{( \, -\infty, 0 \, ]}( f(x,\tilde{z}) ) \, \right]}{
\mathbb{P}_{\tilde{z}}\left( \, f(x,\tilde{z}) \leq 0 \, \right)}
} \, = \, \displaystyle{
\frac{\mathbb{E}_{\tilde{z}}\left[ \, \phi(x,\tilde{z}) \, \onebld_{( \, -\infty, 0 \, ]}( f(x,\tilde{z}) ) \, \right]}{
\mathbb{E}_{\tilde{z}}\left[ \, \onebld_{( \, -\infty, 0 \, ]}( f(x,\tilde{z}) ) \, \right]},
} \]
which again involves the closed Heaviside function.
 {Conditional expectations have applications in modeling low probability--high consequence accidents \cite{Sherali1997}
and in the semi-supervised structured classification problem \cite{ZhengChang2016}}.
The understanding and treatment of the deterministic problem (\ref{eq:original POP}) provides useful insights
for studying stochastic programs with composite indicator functions such as the optimization with chance constraints
involving discontinuous piecewise functions
and conditional expectation functions that lead to fractional expectation functions, where the Heaviside functions are prominently present.




\section{Equivalent Formulations: Global Minimizers} \label{sec:global minimizers}

Before discussing stationarity solutions and local minimizers, we introduce classes of
the problem (\ref{eq:original POP}) within which a (globally) optimal solution exists
and there is an equivalence of the problem with several lifted formulations employing additional variables and special constraints.
Equivalence means that
there is a one-to-one correspondence between the optimal solutions of the problems and their optimal objective values are equal.
The purpose of this section is twofold: (a) to show that the problem (\ref{eq:original POP}) has an optimal solution and is related to two well-studied
problems in the literature under some sign restrictions on the family of functions $\{ \varphi_k,g_k \}_{k=1}^K$ and $\{ \phi_{\ell},h_{\ell} \}_{\ell=1}^L$,
and (b) to support the sign conditions as a reasonable assumption to be imposed in subsequent sections.  Since these
equivalent formulations are for global minimizers whose computation is practically prohibitive if not impossible, and
since they are not the focus in the later sections, we present the result below without the Heaviside constraint;
i.e., for the following problem only:
\begin{equation} \label{eq:original POP no Heaviside constraint}
\displaystyle{
\operatornamewithlimits{\mbox{\bf minimize}}_{x \in X}
} \ \Phi(x) \, \triangleq \, c(x) +
\displaystyle{
\sum_{k=1}^K
} \, \varphi_k(x) \, \onebld_{( \, 0,\infty \, )}( g_k(x)).
\end{equation}
In the proof, we let $( \bullet )_{\min}$ denote the minimum objective value of the referenced problem; $( t )_{\pm} \triangleq \max( \pm t, 0 )$
be the nonnegative and nonpositive part of a scalar $t$; and let $\perp$ denote
the perpendicularity notation, which in the present context means the complementary slackness between the involved expressions.

\begin{proposition} \label{pr:lsc case} \rm
Let $X$ be a compact set.  Suppose that
the functions $c$ 
and $\{ \varphi_k, g_k \}_{k=1}^K$ are continuous.  The following two statements hold:

\gap

(A)  Under the following sign restriction:

\gap
%
%
%
$\bullet $  for every $k \in [ K ]$,
the function $\varphi_k$ is nonnegative on the set $X \cap g_k^{-1}(0)$,

\gap

the problem (\ref{eq:original POP no Heaviside constraint}) has an optimal solution; moreover, it is
equivalent to:

\gap

--- {\bf MPCC-1:}
\begin{equation} \label{eq:equivalent MPCC-1}
\begin{array}{rl}
\displaystyle{
\operatornamewithlimits{\mbox{\bf minimize}}_{x \, \in \, X; \, s}
} & \Phi_{\rm MPCC1}(x,s) \, \triangleq \, c(x) + \displaystyle{
\sum_{k=1}^K
} \, \varphi_k(x) \, s_k \\ [0.2in]
\mbox{\bf subject to} & 0 \, \leq s_k \, \perp \, ( \, g_k(x) \, )_- \, \geq \, 0, \hspace{0.4in} \ \forall \, k \, \in \, [ \, K \, ] \\ [0.1in]
\mbox{\bf and} & 0 \, \leq \, 1 - s_k \, \perp \, ( \, g_k(x) \, )_+ \, \geq \, 0, \epc \ \forall \, k \, \in \, [ \, K \, ].
\end{array}
\end{equation}
(B)  Under the following strengthened sign restriction:

\gap

$\bullet $  for every $k \in [ K ]$,
the function $\varphi_k$ is nonnegative on the set $X \cap g_k^{-1}( \, -\infty, 0 \, ]$,

\gap

the problem (\ref{eq:original POP no Heaviside constraint}) is equivalent to either one of the following two problems:

\gap

--- {\bf MPCC-2:}
\begin{equation} \label{eq:equivalent MPCC}
\begin{array}{rl}
\displaystyle{
\operatornamewithlimits{\mbox{\bf minimize}}_{x \, \in \, X; \, s}
} & \Phi_{\rm MPCC2}(x,s) \, \triangleq \, c(x) + \displaystyle{
\sum_{k=1}^K
} \, \varphi_k(x) \, s_k \\ [0.2in]
\mbox{\bf subject to}
& 0 \, \leq \, 1 - s_k \, \perp \, ( \, g_k(x) \, )_+ \, \geq \, 0, \epc \ \forall \, k \, \in \, [ \, K \, ] \\ [0.1in]
\mbox{\bf and} & s \, \in \, [ \, 0,1 \, ]^K \epc \mbox{(continuous variables)};
\end{array}
\end{equation}
(note the absence of the constraint $s_k ( \, g_k(x) \, )_- = 0$ compared to (\ref{eq:equivalent MPCC-1}));

\gap

--- {\bf on-off constraints:}
\begin{equation} \label{eq:original on-off}
\begin{array}{ll}
\displaystyle{
\operatornamewithlimits{\mbox{\bf minimize}}_{x \ \in \, X; \, z}
} & \Phi_{\rm on/off}(x,z) \, \triangleq \, c(x) + \displaystyle{
\sum_{k=1}^K
} \, \varphi_k(x) \, ( 1 - z_k ) \\ [0.2in]
\mbox{\bf subject to}
& g_k(x) \, \leq \, 0 \epc \mbox{if $z_k = 1$}, \ \forall \, k \, \in \, [ K ] \\ [0.1in]
\mbox{\bf and} & z \, \in \, \{ \, 0,1 \}^K \epc \mbox{(binary variables)}.
\end{array} \end{equation}
\end{proposition}

\begin{proof}   We claim that the objective $\Phi$ is lower semicontinuous on $X$ by showing that each product
$\varphi_k(x) \, \onebld_{( \, 0,\infty \, )}( g_k(x))$
is lower semicontinuous on $X$ under the sign assumption on the element functions
$\{ \varphi_k \}_{k=1}^K$.  For this purpose, it suffices to show if $\varphi$ and $g$ are two continuous functions such that
$\varphi$ is nonnegative on $X \cap g^{-1}(0)$, then the level set:
\[
L(\alpha) \, \triangleq \, \{ \, x \, \in \, X \, \mid \, \varphi(x) \, \onebld_{( \, 0,\infty \, )}(g(x))  \, \leq \, \alpha \, \}
\]
is closed for all scalars $\alpha \in \mathbb{R}$.  Let $\{ x^{\, \nu} \}$ be a given sequence in $L(\alpha)$ converging to a limit $\bar{x} \in X$.
There are 2 cases to consider:

\gap

$\bullet $ $\alpha \geq 0$:  We must have for every $\nu$, either $g(x^{\, \nu}) \leq 0$ or [ $g(x^{\, \nu}) > 0$ and $\varphi(x^{\, \nu}) \leq \alpha$ ].
Then the limit $\bar{x}$ satisfies: either $g(\bar{x}) \leq 0$ or [ $g(\bar{x}) > 0$ and $\varphi(\bar{x}) \leq \alpha$ ].  This shows that
$\bar{x} \in L(\alpha)$, without requiring the sign restriction on $\varphi$.

\gap

$\bullet $ $\alpha < 0$: We must have [ $g(x^{\, \nu}) > 0$ and $\varphi(x^{\, \nu}) \leq \alpha$ ].  By the sign assumption on $\varphi$, the limit $\bar{x}$
must satisfy: $g(\bar{x}) > 0$ and $\varphi(\bar{x}) \leq \alpha$; so $\bar{x} \in {\cal L}(\alpha)$, completing the proof of the closedness of
$L(\alpha)$.

\gap

Consequently, the problem (\ref{eq:original POP no Heaviside constraint}) is a minimization problem of a lower semicontinous function on a compact set; thus it has
an optimal solution.  Omitting the proof of equivalence with the problem (\ref{eq:equivalent MPCC-1}), we directly show
the equivalence of the two problems (\ref{eq:original POP no Heaviside constraint}) and (\ref{eq:equivalent MPCC}).  Let $x \in X$ be arbitrary.
Then the pair $(x,s)$ is
feasible to (\ref{eq:equivalent MPCC}), where 
$s_k \triangleq \onebld_{( \, 0,\infty \, )}(g_k(x))$
for all $k \in [ K ]$.  Hence
 {the minimum objective value of (\ref{eq:original POP no Heaviside constraint}) is no less than that of
(\ref{eq:equivalent MPCC})}.
Conversely,
if $(x,s)$ is feasible to (\ref{eq:equivalent MPCC}), then $\Phi_{\rm MPCC2}(x,s) \geq \Phi(x)$ by the sign restriction of $\varphi_k$
on $X \cap g_k^{-1}( \, -\infty,0 \, ]$.
Therefore, equality holds and an optimal solution
of one problem readily yields an optimal solution of the other.  To show the equivalence of (\ref{eq:original POP no Heaviside constraint})
and (\ref{eq:original on-off}), let $x \in X$ be arbitrary.  Then the pair $(x,z)$ is
feasible to (\ref{eq:original on-off}), where 
$1 - z_k \triangleq \onebld_{( \, 0,\infty \, )}(g_k(x))$
for all $k \in [ K ]$.   Hence $(\ref{eq:original POP no Heaviside constraint})_{\min} \geq (\ref{eq:original on-off})_{\min}$.  Conversely,
if $(x,z)$ is feasible to (\ref{eq:original on-off}), then $\Phi_{\rm on/off}(x,z) \geq \Phi(x)$.   Hence equality holds and there is a one-to-to
correspondence between the optimal solutions of these two problems.
\end{proof}

 {Being instances of a mathematical program
with complementarity constraints (MPCC), the formulations (\ref{eq:equivalent MPCC-1}) and (\ref{eq:equivalent MPCC}) are similar to
the equivalent formulations of an $\ell_0$-minimization problem used in \cite{FMPangSWachter18}.
These two MPCC's (\ref{eq:equivalent MPCC-1}) and (\ref{eq:equivalent MPCC}) are special instances of a mathematical program with vanishing
constraints for which there is an extensive literature; see \cite{AchtzigerKanzow08,DHMigot19,IzmailovSolodov09,HKanzowSchwartz12} which contain
many more references.  With an emphasis toward more general constraints, these references have paid significant attention to constraint
qualifications \cite{FukushimaPang98} and their consequences for stationarity conditions and regularization methods.  In contrast, our study bypasses
such complementarity constraint qualifications and aims to analyze a kind of stationary solutions defined directly on the discontinuous piecewise
problem (\ref{eq:original POP}).  The formulation (\ref{eq:original on-off})
is one with on-off constraints described by indicator variables \cite{BLTrWiese15,HBCornuejols12}.  The method of proof of the
equivalence of the two problems (\ref{eq:equivalent MPCC}) and (\ref{eq:original on-off}) is fairly elementary; the same methodology
is also discussed and used, e.g., in \cite{PAhmedShapiro09,KNNorkin13}.}

\gap

 {The principal purpose of presenting Proposition~\ref{pr:lsc case} is to highlight the two sign assumptions in (A) and (B).
They delimit the applicability of the theory in the rest of the paper.  The sign conditions can be constructively verified (if desired)
when the functions $\varphi_k$ are convex and $g_k$ are affine (for (A)) or convex (for (B)).  In general, these assumptions can be
formulated equivalently as the optimum objective values of the minimization problems being nonnegative:
\[ \displaystyle{
\operatornamewithlimits{\mbox{\bf minimize}}_{x \in X}
} \ \varphi_k(x) \ \mbox{ \bf subject to } \ g_k(x) \, = \, 0 \epc | \epc
\displaystyle{
\operatornamewithlimits{\mbox{\bf minimize}}_{x \in X}
} \ \varphi_k(x) \ \mbox{ \bf subject to } \ g_k(x) \, \leq \, 0.
\]
respectively.  These are convex programs under stated stipulations of $g_k$ and $\varphi_k$.
}

\section{Pseudo Bouligand Stationarity} \label{sec:pseudo}

With the discontinuity of the (open) Heaviside function it is not easy to derive transparent necessary conditions for
a local minimizer of the problem (\ref{eq:original POP}), let alone computing it.  Although there are various one-sided
directional derivatives (such as those of the Dini kind), and more generally, the subderivatives in modern variational
analysis \cite{RockafellarWets09} that one may apply to the objective function $\Phi(x)$, these derivatives
are defined for general functions; in particular, their specializations to functions such as $\Phi(x)$ in (\ref{eq:original POP})
do not immediately yield useful insights about the problem without carefully unwrapping the details of the derivatives.  Instead,
our approach herein is based on elementary one-sided directional derivatives and basic optimization theory. One immediate benefit
of our approach is that it handles constraints with the Heaviside functions at much ease, unlike the variational approach which requires
the constraints to be converted to extended-valued functions embedded in the objective.  Another important point to make is that we
aim to connect the theory with computations; namely, we want to ensure two goals of the defined solution concepts: (i) they are
computationally achievable, and (ii) they are not over-relaxed.
As the background for the subsequent Definition~\ref{df:pseudo concepts}, we provide a result \cite[Proposition~3]{GomezHePang21}
that was stated for a global minimizer in the reference, which we broaden to allow for a local minimizer.

\begin{proposition} \label{pr:Necessary for glomin} \rm
Let $X \subseteq \mathbb{R}^n$, $\theta : {\cal O} \supseteq X \to \mathbb{R}$, and
$\gamma > 0$ be given.  If $\bar{x}$ is a global (local) minimizer of the problem:
\begin{equation} \label{eq:ell0 general}
\displaystyle{
\operatornamewithlimits{\mbox{\bf minimize}}_{x \in X}
} \ \theta(x) + \gamma \, \displaystyle{
\sum_{i=1}^n
} \, | \, x_i \, |_0,
\end{equation}
then $\bar{x}$ is a global (local, respectively) minimizer of the problem:
\begin{equation} \label{eq:restricted zero-problem}
\displaystyle{
\operatornamewithlimits{\mbox{\bf minimize}}_{x \in X}
} \ \theta(x) \epc \mbox{\bf subject to } \ x_{{\cal A}_0(\bar{x})} \, = \, 0,
\end{equation}
where ${\cal A}_0(\bar{x}) \triangleq \left\{ i \in [ n ] \mid \bar{x}_i = 0 \right\}$.
\hfill $\Box$
\end{proposition}

Restricted to problems with $\theta$ convex, a vector $\bar{x}$ that is a global minimizer of
(\ref{eq:restricted zero-problem}) is termed a ``pseudo-minimizer'' of (\ref{eq:ell0 general}) in the reference.
It turns out that such minimizers are the points of attraction of sequences of stationary solutions
of the ``folded concave'' approximations of the $\ell_0$-problem; see \cite[Proposition~5]{GomezHePang21}.
The main idea of the problem (\ref{eq:restricted zero-problem}) to deal with the discontinuous
function $| \bullet |_0$ is to transfer, or ``pull down'' the discontinuity of the objective to a constraint.
We plan to apply this idea to the problem (\ref{eq:original POP}); nevertheless the resulting ``pulled-down''
extension of (\ref{eq:restricted zero-problem}) will have nonconvex objectives and constraints.  This necessitates
us to recall the definition of a Bouligand stationary (B-stationary in short)
solution of a B-differentiable program.   For an abstract optimization problem:  $\displaystyle{
\operatornamewithlimits{\mbox{\bf minimize}}_{x \in \wh{X}}
} \, \theta(x)$, where $\wh{X}$ is a closed set in $\mathbb{R}^n$ and $\theta : {\cal O} \to \mathbb{R}$ is a B-differentiable function,
a vector $\bar{x} \in \wh{X}$ is a Bouligand stationary solution of $\theta$ on $\wh{X}$ \cite[Definition~6.1.1]{CuiPang2021}
if $\theta^{\, \prime}(\bar{x};v) \geq 0$ for all $v \in {\cal T}(\wh{X};\bar{x})$, where
\[
{\cal T}(\wh{X};\bar{x}) \, \triangleq \, \left\{ \, v \, \in \, \mathbb{R}^n \, \left| \, \exists \ \{ x^{\, \nu} \} \subset \wh{X} \mbox{ converging to $\bar{x}$ and
$\{ \tau_{\nu} \} \downarrow 0$ such that $v \, = \, \displaystyle{
\lim_{\nu \to \infty}
} \ \displaystyle{
\frac{x^{\, \nu} - \bar{x}}{\tau_{\nu}}
}$} \, \right. \right\}
\]
is the tangent cone of $\wh{X}$ at $\bar{x}$.  If $\wh{X}$ is additionally convex, B-stationarity reduces to d-stationary (``d'' for directional),
which is: $\theta^{\, \prime}(\bar{x};x - \bar{x}) \geq 0$ for all $x \in \wh{X}$.

\gap

For the problem (\ref{eq:original POP}), we define the following six index
sets corresponding to a given $\bar{x} \in X$:
\[ \begin{array}{ll}
{\cal K}_>(\bar{x}) \, \triangleq \, \{ \, k \, \in \, [ K ] \, \mid \, g_k(\bar{x}) \, > \, 0 \, \}; &
{\cal K}_=(\bar{x}) \, \triangleq \, \{ \, k \, \in \, [ K ] \, \mid \, g_k(\bar{x}) \, = \, 0 \, \} \\ [0.1in]
{\cal K}_<(\bar{x}) \, \triangleq \, \{ \, k \, \in \, [ K ] \, \mid \, g_k(\bar{x}) \, < \, 0 \, \} \\ [0.1in]
{\cal L}_>(\bar{x}) \, \triangleq \, \{ \, \ell \, \in \, [ L ] \, \mid \, h_{\ell}(\bar{x}) \, > \, 0 \, \}; &
{\cal L}_=(\bar{x}) \, \triangleq \, \{ \, \ell \, \in \, [ L ] \, \mid \, h_{\ell}(\bar{x}) \, = \, 0 \, \} \\ [0.1in]
{\cal L}_<(\bar{x}) \, \triangleq \, \{ \, \ell \, \in \, [ L ] \, \mid \, h_{\ell}(\bar{x}) \, < \, 0 \, \}.
\end{array} \]
With these index sets, we define the ``pulled-down'', or ``pseudo stationarity'' problem at $\bar{x}$:
\begin{equation} \label{eq:barx stationary problem}
\begin{array}{l}
\displaystyle{
\operatornamewithlimits{\mbox{\bf minimize}}_x
} \  \Phi(x;\bar{x}) \, \triangleq \, c(x) + \displaystyle{
\sum_{k \, \in \, {\cal K}_>(\bar{x})}
} \, \varphi_k(x) \\ [0.3in]
\left. \begin{array}{ll}
\mbox{\bf subject to} & x \, \in \, X; \ \displaystyle{
\sum_{\ell \, \in \, {\cal L}_>(\bar{x})}
} \, \phi_{\ell}(x) \, \leq \, b \\ [0.2in]
& g_k(x) \, \leq \, 0 \epc \forall \, k \in {\cal K}_=(\bar{x}) \, \cup \, {\cal K}_<(\bar{x}) \triangleq {\cal K}_{\leq}(\bar{x})
\\ [0.1in]
& g_k(x) \, \geq \, 0 \epc \forall \, k \in {\cal K}_>(\bar{x}) \\ [0.1in]
& h_{\ell}(x) \, \leq \, 0 \epc \forall \, \ell \in {\cal L}_=(\bar{x}) \, \cup \, {\cal L}_<(\bar{x}) \triangleq {\cal L}_{\leq}(\bar{x}) \\ [0.1in]
\mbox{\bf and} & h_{\ell}(x) \, \geq \, 0 \epc \forall \, \ell \in {\cal L}_>(\bar{x})
\end{array} \right\} \epc \mbox{\begin{tabular}{l}
feasible set \\ [5pt]
denoted $S_{\rm ps}(\bar{x})$.
\end{tabular}}
\end{array} \end{equation}
Quite different from the three problems in Proposition~\ref{pr:lsc case},
the above problem is obtained by pulling the indicator functions out of the objective function and the functional constraint
and imposing constraints restricting the variable $x$ according to those defined by the functions $\{ g_k,h_{\ell} \}$
and satisfied by the vector $\bar{x}$ on hand.
Specialized to the $\ell_0$-optimization problem (\ref{eq:ell0 general}),
the above pulled-down problem reduces to (\ref{eq:restricted zero-problem}).
 {Clearly, if $\bar{x}$ is a globally optimal solution to (\ref{eq:original POP}) and provided that
$\varphi_k$ is nonnegative on $X \cap g_k^{-1}(0)$ for $k \in {\cal K}_>(\bar{x})$ and $\phi_{\ell}$ is nonnegative on
$X \cap h_{\ell}^{-1}(0)$ for $\ell \in {\cal L}_>(\bar{x})$, then $\bar{x}$ is globally optimal for
(\ref{eq:barx stationary problem}); however, the converse is generally not true because $S_{\rm ps}(\bar{x})$ is at best
only a subset of the feasible set of (\ref{eq:original POP}).  This failed equivalence is not surprising for the main reason that
(\ref{eq:barx stationary problem}) is defined at a given $\bar{x}$; this definition is for the purpose of answering the question:
what is a necessary condition for a given feasible vector of (\ref{eq:original POP}) to be its local
minimizer?}

\gap

The following definition specifies a pseudo solution of
(\ref{eq:original POP}) as a fixed point of the self-defined ``locmin'' or ``B-stationarity'' mapping.  Note that with $\bar{x}$ given,
problem (\ref{eq:barx stationary problem}) is a standard nonlinear program;  {thus, a B-stationary point of the problem
is a vector $\wh{x} \in S_{\rm ps}(\bar{x})$ such that $\Phi(\bullet;\bar{x})^{\, \prime}(\wh{x};v) \geq 0$ for
all $v \in {\cal T}(S_{\rm ps}(\bar{x});\wh{x})$.}

\begin{definition} \label{df:pseudo concepts} \rm
A vector $\bar{x} \in X$ is said to be a

\gap

$\bullet $ {\sl pseudo B-stationary point} of (\ref{eq:original POP}) if $\bar{x}$ 
is a B-stationary point of the pulled-down problem (\ref{eq:barx stationary problem}). 

\gap

$\bullet $ {\sl pseudo local minimizer} of (\ref{eq:original POP}) if $\bar{x}$ 
is a local minimizer of the pulled-down problem (\ref{eq:barx stationary problem}). 
\hfill $\Box$
\end{definition}

 {Besides its dictionary meaning of being not real, the adjective ``pseudo'' carries the hiddent meaning that the
concepts are defined in terms of an auxiliary problem induced by the candidate solution on hand.}
The following result shows that the above pseudo conditions are necessary for a local minimizer of (\ref{eq:original POP}), thus providing a first
step in computing a promising candidate for a local minimum of the problem.

\begin{proposition} \label{pr:pseudo as necessary} \rm
Let $c$, $\{ \varphi_k,g_k \}_{k=1}^K$, and $\{ \phi_{\ell},h_{\ell} \}_{\ell=1}^L$ be B-differentiable functions on
the open set ${\cal O}$ containing the closed set $X$.
Among the following three statements for a vector $\bar{x} \in X$, it holds that (a) $\Rightarrow$ (b) $\Rightarrow$ (c):

\gap

(a) $\bar{x}$ is a local minimizer of (\ref{eq:original POP}).

\gap

(b) $\bar{x}$ is a pseudo local minimizer of (\ref{eq:original POP}).

\gap

(c) $\bar{x}$ is a pseudo B-stationary point of (\ref{eq:original POP}).
\end{proposition}

\begin{proof} It suffices to show that the local minimizer of problem (\ref{eq:original POP}) must be a pseudo local minimizer.
By the continuity of the functions
$\{ g_k \}_{k=1}^K$,
there exists a neighborhood ${\cal N}$ of $\bar{x}$ such that for all $x \in {\cal N}$, it holds that:
$g_k(x) < 0$ for all $k \in {\cal K}_<(\bar{x})$ and $g_k(x) > 0$ for all $k \in {\cal K}_>(\bar{x})$;
and similarly for the $\{ \phi_{\ell},h_{\ell} \}$-functions.
We may restrict the neighborhood ${\cal N}$ so that $\bar{x}$ is a minimizer of $\Phi$ on $X \cap {\cal N}$.
Hence if $x \in {\cal N}$ is feasible to (\ref{eq:barx stationary problem}), then $x$ is also feasible to (\ref{eq:original POP})
because $\displaystyle{
\sum_{\ell=1}^L
} \, \phi_{\ell}(x) \onebld_{(0,\infty)}(h_{\ell}(x)) = \displaystyle{
\sum_{\ell \in {\cal L}_>(\bar{x})}
} \, \phi_{\ell}(x)$. Thus we have
\[
\Phi(x;\bar{x}) \, = \, \Phi(x) \, \geq \, \Phi(\bar{x}) \, = \, \Phi(\bar{x};\bar{x})
\]
showing that $\bar{x}$ is a local minimizer of the problem (\ref{eq:barx stationary problem}); thus (b) holds.
\end{proof}

An important point of Proposition~\ref{pr:pseudo as necessary} is that the two pseudo conditions are necessary
for a local minimizer of (\ref{eq:original POP}).  With the two approaches described in the two later sections
and with the aid of a host of existing algorithms for practically implementing the approaches (as mentioned in the closing
of the Introduction), the computation of a pseudo B-stationary solution can be accomplished
by iterative algorithms for large classes of functions.  This is in contrast to the computation of a pseudo local minimizer
which in general is a daunting, if not impossible task.  It is therefore natural to ask whether there are nonconvex nondifferentiable
classes of problems
for which a pseudo B-stationary point must be pseudo locally minimizing.  It turns out that the answer is affirmative based
on the convex-like property defined as follows;
see \cite[Section~4.1]{CuiLiuPang21}.  A function $f : \mathbb{R}^n \to \mathbb{R}$ is said to be
{\sl convex-like} near a vector $\bar{x}$ if there exists a neighborhood ${\cal N}$
of $\bar{x}$ such that
\[
f(x) \, \geq \, f(\bar{x}) + f^{\, \prime}(\bar{x};x - \bar{x}), \epc
\forall \, x \, \in \, {\cal N}.
\]
A very broad class of convex-like functions consists of the following 3-layer composite functions:
\begin{equation} \label{eq:generic fk}
f(x) \, \triangleq \, \varphi \circ \theta \, \circ \, \psi(x),
\end{equation}
where $\varphi : \mathbb{R} \to \mathbb{R}$ is piecewise affine and nondecreasing;
$\theta : \mathbb{R} \to \mathbb{R}$ is convex, and
$\psi : \mathbb{R}^n \to \mathbb{R}$ is piecewise affine; see \cite[Lemma~10]{CuiLiuPang21}.
In particular, piecewise affine functions and convex functions are convex-like near any point.
Under the convexity-like conditions, we aim to specialize \cite[Proposition~9]{CuiLiuPang21} to
the problem (\ref{eq:barx stationary problem}).  For this purpose, let $\bar{x} \in S_{\rm ps}(\bar{x})$.
We then have the inclusion:
\begin{equation} \label{eq:tangent subset linearized}
{\cal T}(S_{\rm ps}(\bar{x});\bar{x}) \, \subseteq \, \left\{ \,
v \, \in \, {\cal T}(X;\bar{x}) \, \left| \, \begin{array}{l}
\displaystyle{
\sum_{\ell \in {\cal L}_>(\bar{x})}
} \, \phi_{\ell}^{\, \prime}(\bar{x};v) \, ( \, \leq \, 0 ) \\ [0.25in]
g_k^{\, \prime}(\bar{x};v) \, \leq \, 0, \epc \forall \, k \in {\cal K}_=(\bar{x}) \\ [0.1in]
h_{\ell}^{\, \prime}(\bar{x};v) \, \leq \, 0, \epc \forall \, \ell \in {\cal L}_=(\bar{x})
\end{array} \right. \right\} \, \triangleq \, {\cal L}(S_{\rm ps}(\bar{x});\bar{x}),
\end{equation}
where the notation $( \, \leq \, 0 )$ means that this constraint is vacuous if the functional
constraint holds as a strict inequality at $\bar{x}$; i.e., if $\displaystyle{
\sum_{\ell=1}^L
} \, \phi_{\ell}(\bar{x}) \, \onebld_{( \, 0,\infty \, )}(h_{\ell}(\bar{x})) = \displaystyle{
\sum_{\ell \, \in \, {\cal L}_>(\bar{x})}
} \, \phi_{\ell}(\bar{x}) < b$.
The equality of the left- and right-hand cones in (\ref{eq:tangent subset linearized})
is the Abadie constraint qualification (ACQ)
for the set $S_{\rm ps}(\bar{x})$ at the member vector $\bar{x}$.
Sufficient conditions for this CQ to hold are known; in particular, either one of the
following two sets of conditions yields the ACQ (see e.g.\ \cite[Section~4.1]{PangRazaAlvarado16}):

\gap

$\bullet $  (piecewise polyhedrality) all the functions
$\{ g_k \}_{k \in {\cal K}_=(\bar{x})}$ and $\{ h_{\ell} \}_{\ell \in {\cal L}_=(\bar{x})}$,
(and $\{ \phi_{\ell} \}_{\ell \in {\cal L}_>(\bar{x})}$ too if the functional constraint is binding at $\bar{x}$)
are piecewise affine; or

\gap

$\bullet $  (directional Slater) there exists a vector $v \in {\cal T}(X;\bar{x})$ such that
\[ \begin{array}{rl}
\displaystyle{
\sum_{\ell \in {\cal L}_>(\bar{x})}
} \, \phi_{\ell}^{\, \prime}(\bar{x};v) & ( \, < \, 0 \, ) \\ [0.2in]
g_k^{\, \prime}(\bar{x};v) & < \, 0, \epc \forall \, k \in {\cal K}_=(\bar{x}) \\ [0.1in]
h_{\ell}^{\, \prime}(\bar{x};v) & < \, 0, \epc \forall \, \ell \in {\cal L}_=(\bar{x}),
\end{array} \]
where $( < 0 )$ has the same meaning as $( \leq 0 )$ when the functional constraint is not binding at $\bar{x}$.
We have the following result which follows readily from  \cite[Proposition~9]{CuiLiuPang21}.

\begin{proposition} \label{pr:convex-like implies} \rm
Under the blanket assumption of (\ref{eq:original POP}), suppose that
$c$, $\{ g_k \}_{k \in {\cal K}_=(\bar{x})}$, $\{ \varphi_k \}_{k \in {\cal K}_>(\bar{x})}$,
$\{ h_{\ell} \}_{\ell \in {\cal L}_=(\bar{x})}$, and $\{ \phi_{\ell} \}_{\ell \in {\cal L}_>(\bar{x})}$,
are all convex-like near $\bar{x}$, which is a pseudo B-stationary solution  of (\ref{eq:original POP}).
If the ACQ holds for the set $S_{\rm ps}(\bar{x})$ at
$\bar{x}$, then $\bar{x}$ is a pseudo local minimizer of (\ref{eq:original POP}).
\hfill $\Box$
\end{proposition}

In terms of the larger cone ${\cal L}(S_{\rm ps}(\bar{x});\bar{x})$,
it follows that a feasible vector $\bar{x}$ of (\ref{eq:original POP}) is a pseudo B-stationary solution
if the following implication holds:
\begin{equation} \label{eq:ACQ stationarity}
v \, \in \, {\cal L}(S_{\rm ps}(\bar{x});\bar{x}) \ \Rightarrow \ \Phi(\bullet;\bar{x})^{\, \prime}(\bar{x};v) \, \geq \, 0;
\end{equation}
or equivalently, if $0 \in \displaystyle{
\operatornamewithlimits{\mbox{\bf argmin}}_{v \, \in \, {\cal L}(S_{\rm ps}(\bar{x});\bar{x})}
} \ \Phi(\bullet;\bar{x})^{\, \prime}(\bar{x};v)$.  When the functions $c$, $\{ g_k \}_{k \in {\cal K}_=(\bar{x})}$,
$\{ \varphi_k \}_{k \in {\cal K}_>(\bar{x})}$,
$\{ h_{\ell} \}_{\ell \in {\cal L}_=(\bar{x})}$, and $\{ \phi_{\ell} \}_{\ell \in {\cal L}_>(\bar{x})}$ are of the difference-of-convex (dc) kind,
then the latter minimization problem in $v$ (with $\bar{x}$ given) is a dc constrained dc program that has been studied in \cite[Section~6.4]{CuiPang2021}.
Under the ACQ for the set $S_{\rm ps}(\bar{x})$ at $\bar{x}$, the implication (\ref{eq:ACQ stationarity}) is necessary and sufficient for
pseudo B-stationarity.  The upshot of this discussion is that with (or without) the ACQ, a feasible vector can in principle be checked for
pseudo B-stationarity by methods existed in the literature; nevertheless, these methods can not be applied to compute a pseudo B-stationary solution.
This is the task in the rest of the paper.

\gap

The pseudo B-stationarity definition can be phrased in an equivalent way, which when strengthened, yields a sufficient condition of a local minimizer
of the problem (\ref{eq:original POP}).  The latter condition is related to the formulations (\ref{eq:equivalent MPCC-1}), (\ref{eq:equivalent MPCC})
and (\ref{eq:original on-off}) albeit with some obvious differences.  First, the auxiliary variables $s$ and $z$ are pulled outside of these formulations;
thereby eliminating the products in the objective functions of these problems.  More interestingly, the result below highlights the difference
between the necessary conditions (in Proposition~\ref{pr:pseudo as necessary}) and the sufficient conditions for a local minimum.  For the former,
the condition is the {\sl existence} of (binary) ``multipliers''; for the latter, the condition is ``for all'' such multipliers.

\begin{proposition} \label{pr:N&S for local} \rm
Let $c$, $\{ \varphi_k,g_k \}_{k=1}^K$, and $\{ \phi_{\ell},h_{\ell} \}_{\ell=1}^L$ be B-differentiable functions on
the open set ${\cal O}$ containing the closed set $X$.  Let $\bar{x} \in X$ be given.  The following three statements hold.

\gap

(a) A necessary condition for $\bar{x}$ to be a local minimizer of (\ref{eq:original POP}) is that {\sl there exist}
(finite) families of binary multipliers
$\{ \xi_k \}_{k \in {\cal K}_=(\bar{x})} \subset \{ 0, 1 \}^{| {\cal K}_=(\bar{x}) |}$ and
$\{ \mu_{\ell} \}_{\ell \in {\cal L}_=(\bar{x})} \subset \{ 0, 1 \}^{| {\cal K}_=(\bar{x}) |}$
such that $\bar{x}$ is a local minimizer of the problem:
\begin{equation} \label{eq:barx necessary locmin}
\begin{array}{l}
\displaystyle{
\operatornamewithlimits{\mbox{\bf minimize}}_x
} \  \Phi(x;\bar{x}) \, \triangleq \, c(x) + \displaystyle{
\sum_{k \, \in \, {\cal K}_>(\bar{x})}
} \, \varphi_k(x) \\ [0.3in]
\left. \begin{array}{ll}
\mbox{\bf subject to} & x \, \in \, X; \ \displaystyle{
\sum_{\ell \, \in \, {\cal L}_>(\bar{x})}
} \, \phi_{\ell}(x) \, \leq \, b \\ [0.2in]
& ( \, 1 - \xi_k \, ) \, g_k(x) \, \leq \, 0 \epc \forall \, k \in {\cal K}_=(\bar{x}) \\ [0.1in]
& g_k(x) \, \leq \, 0 \hspace{0.2in} \forall \, k \, \in \, {\cal K}_<(\bar{x}) \\ [0.1in]
& g_k(x) \, \geq \, 0 \epc \forall \, k \in {\cal K}_>(\bar{x}) \\ [0.1in]
& ( \, 1 - \mu_{\ell} \, ) \, h_{\ell}(x) \, \leq \, 0 \epc \forall \, \ell \, \in \, {\cal L}_=(\bar{x}) \\ [0.1in]
& h_{\ell}(x) \, \leq \, 0 \epc \forall \, \ell \, \in \, {\cal L}_<(\bar{x}) \\ [0.1in]
\mbox{\bf and} & h_{\ell}(x) \, \geq \, 0 \epc \forall \, \ell \, \in \, {\cal L}_>(\bar{x})
\end{array} \right\} \epc \mbox{\begin{tabular}{l}
 {denoted as $\wh{S}_{ps}(\bar x; \xi, \mu)$}, \\ [5pt]
same as $S_{\rm ps}(\bar{x})$ except for \\ [5pt]
the constraints indexed by \\ [5pt]
${\cal K}_=(\bar{x}) \times {\cal L}_=(\bar{x})$.
\end{tabular}}
\end{array}
\end{equation}

(b) Conversely, suppose that in a neighborhood of $\bar{x}$, the products $\varphi_k \, [ g_k ]_+$ and
$\phi_{\ell} \, [ h_{\ell} ]_+$ are nonnegative for all pairs $( k,\ell )$ in ${\cal K}_=(\bar{x}) \times {\cal L}_=(\bar{x})$.
Then a sufficient condition for $\bar{x}$ to be a local minimizer of (\ref{eq:original POP}) is that {\sl for all}
$\{ \xi_k \}_{k \in {\cal K}_=(\bar{x})} \subset \{ 0, 1 \}^{| {\cal K}_=(\bar{x}) |}$ and
$\{ \mu_{\ell} \}_{\ell \in {\cal L}_=(\bar{x})} \subset \{ 0, 1 \}^{| {\cal L}_=(\bar{x}) |}$,
$\bar{x}$ is a local minimizer of (\ref{eq:barx necessary locmin}).

\gap

(c) Alternatively, suppose that $\varphi_k(\bar{x})$ and $\phi_{\ell}(\bar{x})$ are nonnegative for
all $( k,\ell )$ in ${\cal K}_=(\bar{x}) \times {\cal L}_=(\bar{x})$.
Then a sufficient condition for $\bar{x}$ to be a local minimizer of (\ref{eq:original POP}) is that {\sl for all}
$\{ \xi_k \}_{k \in {\cal K}_=(\bar{x})} \subset \{ 0, 1 \}^{| {\cal K}_=(\bar{x}) |}$ and
$\{ \mu_{\ell} \}_{\ell \in {\cal L}_=(\bar{x})} \subset \{ 0, 1 \}^{| {\cal L}_=(\bar{x}) |}$,
$\bar{x}$ is a local minimizer of
\begin{equation} \label{eq:barx sufficient locmin}
\begin{array}{l}
\displaystyle{
\operatornamewithlimits{\mbox{\bf minimize}}_x
} \  \Phi_{\geq}^{\xi}(x;\bar{x}) \, \triangleq \, c(x) + \displaystyle{
\sum_{k \, \in \, {\cal K}_>(\bar{x})}
} \, \varphi_k(x) + \underbrace{\displaystyle{
\sum_{k \, \in \, {\cal K}_=(\bar{x})}
} \, \xi_k \, \varphi_k(x)}_{\mbox{extra term with multiplier}} \\ [0.5in]
\left. \begin{array}{ll}
\mbox{\bf subject to} & x \, \in \, X; \ \displaystyle{
\sum_{\ell \, \in \, {\cal L}_>(\bar{x})}
} \, \phi_{\ell}(x) + \underbrace{\displaystyle{
\sum_{\ell \, \in \, {\cal L}_=(\bar{x})}
} \, \mu_{\ell} \, \phi_{\ell}(x)}_{\mbox{extra term with multiplier}} \, \leq \, b \\ [0.5in]
& ( \, 1 - \xi_k \, ) \, g_k(x) \, \leq \, 0 \epc \forall \, k \in {\cal K}_=(\bar{x}) \\ [0.1in]
& g_k(x) \, \leq \, 0 \hspace{0.2in} \forall \, k \, \in \, {\cal K}_<(\bar{x}) \\ [0.1in]
& g_k(x) \, \geq \, 0 \epc \forall \, k \in {\cal K}_>(\bar{x}) \\ [0.1in]
& ( \, 1 - \mu_{\ell} \, ) \, h_{\ell}(x) \, \leq \, 0 \epc \forall \, \ell \, \in \, {\cal L}_=(\bar{x}) \\ [0.1in]
& h_{\ell}(x) \, \leq \, 0 \epc \forall \, \ell \, \in \, {\cal L}_<(\bar{x}) \\ [0.1in]
\mbox{\bf and} & h_{\ell}(x) \, \geq \, 0 \epc \forall \, \ell \, \in \, {\cal L}_>(\bar{x})
\end{array} \right\} .
\end{array}
\end{equation}
\end{proposition}

\begin{proof} The necessary condition in (a) is clear because we can let $\xi_k = 0 = \mu_{\ell}$ for all
$( k,\ell ) \in {\cal K}_=(\bar{x}) \times {\cal L}_=(\bar{x})$.  For the sufficiency, suppose that $\bar{x}$ is a local
minimizer of (\ref{eq:barx necessary locmin}) for all tuples $\{ \xi_k \}_{k \in {\cal K}_=(\bar{x})}$
and $\{ \mu_{\ell} \}_{\ell \in {\cal L}_=(\bar{x})}$ as stated.
 {We first note that $\bar{x}$ must be feasible to (\ref{eq:original POP}) because
\[
\displaystyle{
\sum_{\ell=1}^L
} \, \phi_{\ell}(\bar{x}) \, \onebld_{( \, 0,\infty \, )}( h_{\ell}(\bar{x})) = \displaystyle{
\sum_{\ell \, \in \, {\cal L}_>(\bar{x})}
} \, \phi_{\ell}(\bar{x}).
\]
}

Since the family of problems (\ref{eq:barx necessary locmin}) is finite, there exists a neighborhood
${\cal N}$ of $\bar{x}$ such that
for all $\{ \xi_k \}_{k \in {\cal K}_=(\bar{x})} \subset \{ 0, 1 \}^{| {\cal K}_=(\bar{x}) |}$ and
$\{ \mu_{\ell} \}_{\ell \in {\cal L}_=(\bar{x})} \subset \{ 0, 1 \}^{| {\cal L}_=(\bar{x}) |}$, if $x \in {\cal N}$ is feasible
to (\ref{eq:barx necessary locmin}), then $\Phi(x;\bar{x}) \geq \Phi(\bar{x};\bar{x}) = \Phi(\bar{x})$.
Without loss of generality, we may assume that this neighborhood is such that for all $x \in {\cal N}$, it holds that
$\varphi_k(x) \, [ g_k(x) ]_+ \geq 0$ and $\phi_{\ell}(x) \, [ h_{\ell}(x) ]_+ \geq 0$ for all
$( k,\ell ) \in {\cal K}_=(\bar{x}) \times {\cal L}_=(\bar{x})$;
$g_k(x) < 0$ for all $k \in {\cal K}_<(\bar{x})$ and $g_k(x) > 0$ for all $k \in {\cal K}_>(\bar{x})$;
and similarly for the $\{ \phi_{\ell},h_{\ell} \}$-functions.  Let
$x \in {\cal N}$ be feasible to (\ref{eq:original POP}).  For each pair $( k,\ell ) \in {\cal K}_=(\bar{x}) \times {\cal L}_=(\bar{x})$,
let $\xi_k \triangleq \onebld_{( \, 0,\infty \, )}(g_k(x))$ and $\mu_{\ell} \triangleq \onebld_{( \, 0,\infty \, )}(h_{\ell}(x))$.
To show that $x$ is feasible to (\ref{eq:barx necessary locmin}) for
this pair $( \xi,\mu )$, it suffices to verify the functional constraint and $( \, 1 - \xi_k \, ) \, g_k(x) \leq 0$ and
$( \, 1 - \mu_{\ell} \, ) \, h_{\ell}(x) \geq 0$ for all
$( k,\ell ) \in {\cal K}_=(\bar{x}) \times {\cal L}_=(\bar{x})$.  The latter is clear.  For the former, we have
\[ \begin{array}{lll}
b & \geq & \displaystyle{
\sum_{\ell=1}^L
} \, \phi_{\ell}(x) \, \onebld_{( \, 0,\infty \, )}(h_{\ell}(x)) \, = \, \displaystyle{
\sum_{\ell \in {\cal L}_>(\bar{x})}
} \, \phi_{\ell}(x) \, \onebld_{( \, 0,\infty \, )}(h_{\ell}(x)) +  \displaystyle{
\sum_{\ell \in {\cal L}_=(\bar{x})}
} \, \phi_{\ell}(x) \, \onebld_{( \, 0,\infty \, )}(h_{\ell}(x)) \\ [0.2in]
& \geq & \displaystyle{
\sum_{\ell \in {\cal L}_>(\bar{x})}
} \, \phi_{\ell}(x) \, \onebld_{( \, 0,\infty \, )}(h_{\ell}(x)) \epc \mbox{because $\phi_{\ell}(x) \, [ h_{\ell}(x) ]_+ \geq 0$ for
$\ell \in {\cal L}_=(\bar{x})$} \\ [0.2in]
& = & \displaystyle{
\sum_{\ell \in {\cal L}_>(\bar{x})}
} \, \phi_{\ell}(x).
\end{array}
\]
Hence it follows that
\[ \begin{array}{lll}
\Phi(x) & = &  {c(x)} + \displaystyle{
\sum_{k \in {\cal K}_>(\bar{x})}
} \, \varphi_k(x) \, \onebld_{( \, 0,\infty \, )}(g_k(x)) +  \displaystyle{
\sum_{k \in {\cal K}_=(\bar{x})}
} \, \varphi_k(x) \, \onebld_{( \, 0,\infty \, )}(g_k(x)) \\ [0.2in]
& \geq & \Phi(x,\bar{x})  \epc \mbox{because $\varphi_k(x) [ g_k(x) ]_+ \geq 0$ for
$k \in {\cal K}_=(\bar{x})$} \\ [0.1in]
& \geq & \Phi(\bar{x};\bar{x}) \, = \, \Phi(\bar{x}),
\end{array} \]
as desired.  Finally, the proof of statement (c) is very similar to that of (b).  Omitting the details,
 {we
simply note that $\bar{x}$, being feasible to (\ref{eq:barx sufficient locmin}), must be feasible to
(\ref{eq:barx necessary locmin}), and thus to (\ref{eq:original POP}).}
\end{proof}

{\bf Remarks.}   {To more precisely connect the pull-down problem (\ref{eq:barx stationary problem})
with the problem (\ref{eq:barx necessary locmin}) for various choices of the binary pairs $( \xi,\mu )$, we note
that the given $\bar{x}$ is a local minimizer of the former if and only if $\bar{x}$ is a local minimizer of the
latter for {\sl some} $(\xi,\mu)$.  In view of statement (b) in Proposition~\ref{pr:N&S for local}, which involves
{\sl for all} $(\xi,\mu)$, it is natural to ask the following question: suppose $\bar{x}$ is a local minimizer
of (\ref{eq:barx stationary problem}), are there sufficient conditions that will ensure $\bar{x}$ to be a local
minimizer of (\ref{eq:original POP})?  Since local optimality involves a neighborhood of $\bar{x}$ wherein
the signs of $\varphi_k(x)[g_k(x)]_+$ and $\phi_{\ell}(x)[ h_{\ell}(x) ]_+$ are relevant, an answer to the question
would essentially reduce to statement (b) in Proposition~\ref{pr:N&S for local}.}

\gap

Incidentally, the difference between parts (b) and (c) in Proposition~\ref{pr:N&S for local} is in their respective assumptions
on the functions $( \varphi_k,g_k )$ and $( \phi_{\ell},h_{\ell} )$ for $( k,\ell ) \in {\cal K}_=(\bar{x}) \times {\cal L}_=(\bar{x})$
and the two resulting problems (\ref{eq:barx necessary locmin}) and (\ref{eq:barx sufficient locmin}).  \hfill $\Box$

\section{The Epigraphical Approach} \label{sec:epi}

Guided by the sign assumption in case (A) of Proposition~\ref{pr:lsc case}, which has two parts for the
problem (\ref{eq:original POP}),

\gap

$\bullet $ for every $k \in [ K ]$, each $\varphi_k$ is nonnegative on the set $X \cap g_k^{-1}(0)$; and

\gap

$\bullet $ for every $\ell \in [ L ]$, each $\phi_{\ell}$ is nonnegative on the set $X \cap h_{\ell}^{-1}(0)$;

\gap

we present
in this section a constructive approach for computing a pseudo B-stationary point of
(\ref{eq:original POP}).  The approach is based on an epigraphical formulation of the problem that lifts it to a higher
dimension.   We recall that the epigraph of a function $f$ on the closed set $S \subseteq {\cal O}$ is the set
\[
{\rm epi}(f;S) \, \triangleq \, \left\{ \, ( \, t,x \, ) \, \in \mathbb{R} \times \, S \, \mid \, t \, \geq \, f(x) \, \right\}.
\]
In what follows, we describe the tangent cone of the epigraph
of a discontinuous product function $\psi \onebld_{( \, 0,\infty \, )}(f)$; the description not only provides insights for the demonstration
of the epigraphical approach to successfully accomplish the computational goal, but is also of independent interest as it relates
to some existing results in the literature.
A basic result of this kind for a locally Lipschitz function is available from \cite[Theorem~2.4.9 part (a)]{Clarke83}; but it is
not applicable to the discontinuous Heaviside function.
An advanced result for an arbitrary function, and for the indicator function in particular,
can be found in \cite[Theorem~8.2]{RockafellarWets09} which is based on the notion of
subderivatives.  Rather than going through the calculation of the latter derivatives for the product function of interest, we present an elementary
derivation that exposes the epigraph of $\psi \onebld_{( \, 0,\infty \, )}(f)$ as the union of two closed sets and highlights the consequence
of the sign condition of $\psi$ on the zero set of $f$.

\begin{proposition} \label{pr:epigraph product} \rm
Let $\psi$ and $f$ be B-differentiable functions on the open set ${\cal O}$ that contains a closed set $S$. Suppose that
$\psi$ is nonnegative on $S \cap f^{-1}(0)$.  Let $\pi(x) \triangleq \psi(x) \, \onebld_{( \, 0,\infty \, )}(f(x))$.  Then,
\begin{equation} \label{eq:epi representation}
{\rm epi}(\pi;S) \, = \,
\underbrace{\left\{ \, ( \, t,x \, ) \, \in \, \mathbb{R} \times S \, \mid \, t \, \geq \, \psi(x), \ f(x) \, \geq \, 0 \, \right\}}_{
\mbox{denoted $E_1$}} \, \cup \, \underbrace{\left( \, \mathbb{R}_+ \times ( \, S \cap f^{-1}( \, -\infty,0 \, ] \, ) \, \right)}_{\mbox{denoted $E_2$}}.
\end{equation}
Thus for any pair $(t,x) \in {\rm epi}(\pi;S)$,
\begin{equation} \label{eq:tangent cone union}
{\cal T}({\rm epi}(\pi;S);(t,x)) \, = \, {\cal T}(E_1;(t,x)) \, \cup \, {\cal T}(E_2;(t,x)),
\end{equation}
where ${\cal T}(E_i;(t,x)) \triangleq \emptyset$ if $( t,x ) \not\in E_i$ for $i = 1,2$.  Moreover, the following three
statements (a), (b), and (c) hold:

\gap

(a) $( t,x ) \in {\rm epi}(\pi;S)$ if and only if $x \in S$ and
\begin{equation} \label{eq:minmax epi}
\underbrace{\min\left( \, \max( \, \psi(x) - t, -f(x) \, ), \ \max( \, f(x), -t \, ) \, \right)}_{\mbox{dc in $(x,t)$ if $\psi$ and $f$ are dc}}
\, \leq \, 0;
\end{equation}
The function on the left side is piecewise affine in $(x,t)$ if $\psi$ and $f$ are piecewise affine.

\gap

(b) For a pair $( \bar{t},\bar{x} ) \in {\rm epi}(\pi;S)$ with $\bar{t} = \pi(\bar{x})$, it holds that

\gap

$\bullet $ if $f(\bar{x}) > 0$, then ${\cal T}({\rm epi}(\pi;S);(\bar{t},\bar{x}) ) =
\left\{ \, ( dt, v ) \, \in \, \mathbb{R} \times {\cal T}(S;\bar{x}) \, \mid \, dt \, \geq \, \psi^{\, \prime}(\bar{x};v) \, \right\}$;

\gap

$\bullet $ if $f(\bar{x}) < 0$, then ${\cal T}({\rm epi}(\pi;S);(\bar{t},\bar{x}) ) = \mathbb{R}_+ \times {\cal T}(S;\bar{x})$;

\gap

$\bullet $ if $f(\bar{x}) = 0 < \psi(\bar{x})$, then
\begin{equation} \label{eq:tangent in E2}
{\cal T}({\rm epi}(\pi;S);(\bar{t},\bar{x}) ) \, \subseteq \,
\left\{ \, ( \, dt,v \, ) \, \in \, \mathbb{R}_+ \times {\cal T}(S;\bar{x}) \, \mid \, f^{\, \prime}(\bar{x};v) \, \leq \, 0 \, \right\}
\end{equation}
with equality holding if the set $\{ x \in S \mid f(x) \leq 0 \}$ satisfies the ACQ at $\bar{x}$;

\gap

$\bullet $ if $f(\bar{x}) = 0 = \psi(\bar{x})$, then
\begin{equation} \label{eq:tangent cone inclusion}
\begin{array}{lll}
{\cal T}({\rm epi}(\pi;S);(\bar{t},\bar{x}) ) & \subseteq & \left\{ \, ( \, dt,v \, ) \, \in \, \mathbb{R} \times {\cal T}(S;\bar{x})
\, \mid \, dt \, \geq \, \psi^{\, \prime}(\bar{x};v), \ f^{\, \prime}(\bar{x};v) \, \geq \, 0 \, \right\} \ \bigcup \\ [0.1in]
& & \left\{ \, ( \, dt,v \, ) \, \in \, \mathbb{R}_+ \times {\cal T}(S;\bar{x}) \, \mid \, f^{\, \prime}(\bar{x};v) \, \leq \, 0 \, \right\};
\end{array} \end{equation}
moreover, if $\psi^{\, \prime}(\bar{x};\bullet)$ is nonnegative on ${\cal T}(S;\bar{x}) \cap ( f^{\, \prime}(\bar{x};\bullet) )^{-1}(0)$,
then the right-hand union in (\ref{eq:tangent cone inclusion}) is equal to
${\rm epi}\left(  \psi^{\, \prime}(\bar{x};\bullet) \, \onebld_{( \, 0,\infty \ )}( f^{\, \prime}(\bar{x};\bullet) ); {\cal T}(S;\bar{x}) \right)$; lastly, if
the two sets $\bar{S}_+ \triangleq \left\{ \, x \in S \mid f(x) \geq 0 \, \right\}$ and $\bar{S}_- \triangleq \left\{ \, x \in S \mid f(x) \leq 0 \, \right\}$
satisfy the ACQ at $\bar{x}$, then
\begin{equation} \label{eq:final tangent of epi}
{\cal T}({\rm epi}(\pi;S);(\bar{t},\bar{x}) ) \, = \,
{\rm epi}\left(  \psi^{\, \prime}(\bar{x};\bullet) \, \onebld_{( \, 0,\infty \ )}( f^{\, \prime}(\bar{x};\bullet) ); {\cal T}(S;\bar{x}) \right).
\end{equation}
(c) For a pair $( \bar{t}, \bar{x} ) \in {\rm epi}(\pi;S)$ with $\bar{t} > \pi(\bar{x})$, it holds that

\gap

$\bullet $ if $f(\bar{x}) \neq 0$ or $f(\bar{x}) = 0 = \psi(\bar{x})$, then
${\cal T}({\rm epi}(\pi;S);(\bar{t},\bar{x}) ) = \mathbb{R} \times {\cal T}(S;\bar{x})$;

\gap

$\bullet $ if $f(\bar{x}) = 0 < \psi(\bar{x}) \neq \bar{t}$, then
\[
{\cal T}({\rm epi}(\pi;S);(\bar{t},\bar{x}) ) \, \left\{ \begin{array}{ll}
= \, \mathbb{R} \times {\cal T}(S;\bar{x}) & \mbox{if $\bar{t} > \psi(\bar{x})$} \\ [0.1in]
\subseteq \, \mathbb{R} \times ( {\cal T}(S;\bar{x}) \cap ( f^{\, \prime}(\bar{x};\bullet) )^{-1}( \, -\infty,0 \, ] )
& \mbox{if $\bar{t} < \psi(\bar{x})$};
\end{array} \right.
\]
$\bullet $ if $f(\bar{x}) = 0 < \psi(\bar{x}) = \bar{t}$, then
\[
\begin{array}{lll}
{\cal T}({\rm epi}(\pi;S);(\bar{t},\bar{x}) ) & \subseteq & \left\{ \, ( \, dt,v \, ) \, \in \, \mathbb{R} \times {\cal T}(S;\bar{x})
\, \mid \, dt \, \geq \, \psi^{\, \prime}(\bar{x};v), \ f^{\, \prime}(\bar{x};v) \, \geq \, 0 \, \right\} \ \bigcup \\ [0.1in]
& & \left\{ \, ( \, dt,v \, ) \, \in \, \mathbb{R} \times {\cal T}(S;\bar{x}) \, \mid \, f^{\, \prime}(\bar{x};v) \, \leq \, 0 \, \right\}.
\end{array} \]
\end{proposition}

\begin{proof}  By definition of the epigraph, we have
\[
{\rm epi}(\pi;S) \, = \,
\left\{ \, ( \, t,x \, ) \, \in \, \mathbb{R} \times S \, \mid \, t \, \geq \, \psi(x), \ f(x) \, > \, 0 \, \right\}
\, \cup \, \left( \, \mathbb{R}_+ \times ( \, S \cap f^{-1}( \, -\infty,0 \, ] \, ) \, \right);
\]
thus ${\rm epi}(\pi;S) \subseteq E_1 \cup E_2$.  Conversely, let $( t,x ) \in E_1$ be such that $f(x) = 0$.  Then $t \geq \psi(x) \geq 0$;
hence $(t,x) \in E_2 \subseteq {\rm epi}(\pi;S)$.  Thus, (\ref{eq:epi representation}) holds.  With the definition that
${\cal T}(E_i;(t,x)) \triangleq \emptyset$ if $( t,x ) \not\in E_i$ for $i = 1,2$, the equality (\ref{eq:tangent cone union}) is clear.
Statement (a) holds because of the following equivalence:
\[ \begin{array}{l}
\left[ \, t \, \geq \, \psi(x) \ \mbox{ and } \ f(x) \, \geq \, 0 \, \right] \ \mbox{ or } \ \left[ \, t \, \geq \, 0 \ \mbox{ and } \ f(x) \, \leq \, 0 \, \right]
\\ [0.1in]
\epc \Leftrightarrow \ [ \, \max( \, \psi(x) - t, -f(x) ) \, \leq \, 0 \, ] \ \mbox{ or } \ [ \, \max( \, f(x), -t \, ) \, \leq \, 0 \ ] \\ [0.1in]
\epc \Leftrightarrow \ \min\left( \, \max( \, \psi(x) - t, -f(x) \, ), \ \max( \, f(x), -t \, ) \, \right) \, \leq \, 0.
\end{array} \]
To prove (b), suppose $f(\bar{x}) > 0$.  Then $f(x) > 0$ for all $x$ sufficiently close to $\bar{x}$.  In this case, the equality
${\cal T}({\rm epi}(\pi;S);(\bar{t},\bar{x}) ) =
\left\{ \, ( dt, v ) \, \in \, \mathbb{R} \times {\cal T}(S;\bar{x}) \, \mid \, dt \, \geq \, \psi^{\, \prime}(\bar{x};v) \, \right\}$
can be proved as follows.  The inclusion $\subseteq$ is straightforward to prove.  For the ``$\supseteq$'' inclusion, let
$( dt, v )$ be such that $v \in {\cal T}(S;\bar{x})$ and $dt \geq \psi^{\, \prime}(\bar{x};v)$.  Let
$v = \displaystyle{
\lim_{\nu \to \infty}
} \,  \displaystyle{
\frac{x^{\, \nu} - \bar{x}}{\tau_{\nu}}
}$ for some sequence $\{ x^{\, \nu} \} \subset S$ converging to $\bar{x}$ and some sequence of scalars $\{ \tau_{\nu} \} \downarrow 0$.
Suppose $dt > \psi^{\, \prime}(\bar{x};v)$.  Then $\bar{t} + \tau_{\nu} dt > \psi(x^{\, \nu})$ for all $\nu$ sufficiently large.  Thus
$( \bar{t} + \tau_{\nu} dt,x^{\nu} ) \in {\rm epi}(\pi;S)$ for all such $\nu$.  This shows that $( dt, v ) \in {\cal T}({\rm epi}(\pi;S);(\bar{t},\bar{x}))$
in this case.  Suppose that $dt = \psi^{\, \prime}(\bar{x};v)$.  Then for every $\varepsilon > 0$,
$( dt + \varepsilon, v ) \in {\cal T}({\rm epi}(\pi;S);(\bar{t},\bar{x}) )$.  Since the tangent cone is closed, we deduce
$( dt, v ) \in {\cal T}({\rm epi}(\pi;S);(\bar{t},\bar{x}))$,
completing the proof of the claimed equality of the two cones.

\gap

The second case where $f(\bar{x}) < 0$ can be similarly argued.  For the third case where $f(\bar{x}) = 0 < \psi(\bar{x})$,
we have $\bar{t} = 0$ and $( \bar{t},\bar{x} ) \in E_2 \setminus E_1$.  So the desired inclusion (\ref{eq:tangent in E2}) and the
equality under the stated CQ both hold easily.

\gap

Consider the last case where $f(\bar{x}) = 0 = \psi(\bar{x})$ so that $( \bar{t},\bar{x}) \in E_1 \cap E_2$.
It is easy to prove that ${\cal T}(E_i;(\bar{t},\bar{x}))$ for $i = 1, 2$
is contained in the two sets on the right-hand side
of (\ref{eq:tangent cone inclusion}), respectively.  If
$\psi^{\, \prime}(\bar{x};\bullet)$ is nonnegative on ${\cal T}(S;\bar{x}) \cap ( f^{\, \prime}(\bar{x};\bullet) )^{-1}(0)$,
the claim that the right-hand union in (\ref{eq:tangent cone inclusion}) is equal to
${\rm epi}\left(  \psi^{\, \prime}(\bar{x};\bullet) \, \onebld_{( \, 0,\infty \ )}( f^{\, \prime}(\bar{x};\bullet); {\cal T}(S;\bar{x}) ) \right)$
can be proved similarly to that of the equality (\ref{eq:epi representation}).  Finally, under the two additional ACQs, the equality
(\ref{eq:final tangent of epi}) holds by the respective representation of ${\cal T}(E_i;(\bar{t},\bar{x}))$.
The proof of statement (c) is similar and omitted.
\end{proof}

Summarizing the various cases in the above proposition, we conclude that there are three ``constraint qualifications''
(in place if needed)
to be satisfied by a pair $( \bar{t},\bar{x} ) \in {\rm epi}(\pi;S)$ in order for the tangent cone
${\cal T}({\rm epi}(\pi;S);(\bar{t},\bar{x}) )$ to have an exact representation:

\gap

(i) $\psi^{\, \prime}(\bar{x};\bullet)$ is nonnegative on ${\cal T}(S;\bar{x}) \cap ( f^{\, \prime}(\bar{x};\bullet) )^{-1}(0)$;

\gap

(ii) the set  $\bar{S}_+ \triangleq \left\{ \, x \in S \mid f(x) \geq 0 \, \right\}$ satisfies the ACQ at $\bar{x}$; and

\gap

(iii) the set $\bar{S}_- \triangleq \left\{ \, x \in S \mid f(x) \leq 0 \, \right\}$
satisfies the ACQ at $\bar{x}$.

\gap

Moreover, refining the above analysis, it is easy to derive sufficient conditions for the set
\[
{\rm epi}(\pi;S) \, = \, \left\{ \, ( t,x ) \, \in \, \mathbb{R} \times S \, \mid \, (\ref{eq:minmax epi}) \mbox{ holds} \, \right\}
\]
to satisfy the ACQ at the pair $( \pi(\bar{x}),\bar{x} )$.

\subsection{The penalized epigraphical formulation}

In addition to the challenges associated with the Heaviside functions, the functional constraint also complicates the analysis and
the design of computational
algorithms for the problem (\ref{eq:original POP}).  There are two ways to address this constraint: one is a direct treatment
as a hard constraint; the other is by a soft penalty with the goal of recovering the satisfaction of the constraint and achieving the stationarity
of the problem.  In this section and Section~\ref{sec:approximation computation}, we adopt the latter penalty approach as it offers a
unified treatment with the composite Heaviside functions occurring in the objective only.  Throughout this subsection,
the functions $\varphi_k$ and $\phi_{\ell}$ are nonnegative on $X \cap g_k^{-1}(0)$ and $X \cap h_{\ell}^{-1}(0)$,
respectively.

\gap

As a first step toward a computationally tractable formulation of (\ref{eq:original POP}), we make the substitutions
\begin{equation} \label{eq:variable substitution}
t_k  \, = \, \varphi_k(x) \onebld_{( \, 0,\infty \, )}(g_k(x)) \ \mbox{ and }  s_{\ell} = \phi_{\ell}(x) \onebld_{( \, 0,\infty \, )}(h_{\ell}(x)),
\epc \forall \, ( k,\ell ) \in [ K ] \times [ L ],
\end{equation}
and relax these definitional equalities to inequalities to obtain the epigraphical constraints:
\[
( t_k,x ) \, \in \, {\rm epi}( \pi_k^{\varphi};X) \epc \mbox{and} \epc ( s_{\ell},x ) \in {\rm epi}( \pi_{\ell}^{\phi};X ),
\]
where
\[
\pi_k^{\varphi}(x) \, \triangleq \, \varphi_k(x) \, \onebld_{( \, 0,\infty \, )}(g_k(x)) \ \mbox{ and } \
\pi_{\ell}^{\phi}(x) \, \triangleq \, \phi_{\ell}(x) \, \onebld_{( \, 0,\infty \, )}(h_{\ell}(x)).
\]
We also relax the functional constraint by penalizing it in the objective.  These maneuvers lead to the following
penalty problem defined for a given parameter $\lambda > 0$: 
\begin{equation} \label{eq:original POP epi form}
\begin{array}{ll}
\displaystyle{
\operatornamewithlimits{\mbox{\bf minimize}}_{x \, \in \, X; \, t; \, s}
} & \Phi_{\lambda}(x,t,s) \triangleq \underbrace{c(x) +
\displaystyle{
\sum_{k=1}^K
} \, t_k}_{\mbox{$\Phi(x)$ in epi-form}} +
\lambda \, \underbrace{\max\left( \, \displaystyle{
\sum_{\ell=1}^L
} \, s_{\ell} - b, \ 0 \, \right)}_{\mbox{\begin{tabular}{ll}
constraint residual \\ [3pt]
fnc.\ in epi-form
\end{tabular}}} \\ [0.35in]
\mbox{\bf subject to} & ( \, t_k,x \, ) \, \in \, {\rm epi}( \pi_k^{\varphi};X ) \ \mbox{ and } \
( \, s_{\ell},x \, ) \, \in \, {\rm epi}( \pi_{\ell}^{\phi};X ), \epc \forall \, ( k,\ell ) \, \in \, [ K ] \times [ L ].
\end{array} \end{equation}
We remark that the $s$-variables and the corresponding
constraints $( s_{\ell},x ) \in {\rm epi}( \pi_{\ell}^{\phi};X )$ are not needed for the problem (\ref{eq:original POP no Heaviside constraint})
where the functional constraint is absent.

\gap

By part (a) of Proposition~\ref{pr:epigraph product}, problem (\ref{eq:original POP epi form}) is equivalent to the following problem
with the epigraphical constraints exposed:
\begin{equation} \label{eq:original POP exposed epi form}
\begin{array}{ll}
\displaystyle{
\operatornamewithlimits{\mbox{\bf minimize}}_{x \, \in \, X; \, t; \, s}
} & \Phi_{\lambda}(x,t,s) \, \triangleq \, c(x) +
\displaystyle{
\sum_{k=1}^K
} \, t_k + \lambda \, \max\left( \, \displaystyle{
\sum_{\ell=1}^L
} \, s_{\ell} - b, \ 0 \, \right) \\ [0.2in]
\mbox{\bf subject to} & \min\left( \, \max( \, \varphi_k(x) - t_k, -g_k(x) \, ), \ \max( \, g_k(x), -t_k \, ) \, \right) \, \leq \, 0,
\epc \forall \, k \, \in \, [ K ] \\ [0.1in]
\mbox{\bf and} & \min\left( \, \max( \, \phi_{\ell}(x) - s_{\ell}, -h_{\ell}(x) \, ), \ \max( \, h_{\ell}(x), -s_{\ell} \, ) \, \right) \, \leq \, 0,
\epc \ \forall \, \ell \, \in \, [ L ].
\end{array} \end{equation}
The latter formulation is the computational workhorse to derive a B-stationary solution of (\ref{eq:original POP epi form}) via
the epigraphical approach.  In what follows,
we show that under a suitable directional conditions on the functions $\phi_{\ell}$ in the functional constraint, a finite value of $\lambda$ exists
such that a B-stationary solution of (\ref{eq:original POP epi form}) is a B-stationary solution of (\ref{eq:barx stationary problem}), hence
a pseudo B-stationary solution of (\ref{eq:original POP}).  The demonstration is carried out in 2 steps, with the first step being the recovery
of the equalities of the $t$-variables in (\ref{eq:variable substitution}).   This step is accomplished by the following lemma.

\begin{lemma} \label{lm:stationary of bivariate} \rm
Let $\theta : {\cal O}_x \times {\cal O}_y \to \mathbb{R}$ be a B-differentiable bivariate variable on the open set ${\cal O}_x \times {\cal O}_y$ containing
the closed set $\wh{X} \subseteq \mathbb{R}^{n+m}$.  If $( \bar{x},\bar{y} )$ is a B-stationary point of $\theta$ on $\wh{X}$, then $\bar{x}$ is a
B-stationary point of $\theta(\bullet,\bar{y})$ on $\wh{X}(\bar{y}) \triangleq \{ \, x \, \mid \, ( x,\bar{y} ) \in \wh{X} \, \}$.  A similar statement holds
for the $y$-variable.
\end{lemma}

\begin{proof}  Indeed, let $v \in {\cal T}(\wh{X}(\bar{y});\bar{x})$ be given such that $v = \displaystyle{
\lim_{\nu \to \infty}
} \, \displaystyle{
\frac{x^{\, \nu} - \bar{x}}{\tau_{\nu}}
}$ for some sequence $\{ x^{\, \nu} \}$ in $\wh{X}(\bar{y})$ converging to $\bar{x}$ and some sequence $\{ \tau_{\nu} \} \downarrow 0$.  It then follows that
$( v,0 ) \in {\cal T}(\wh{X};(\bar{x},\bar{y}))$ and
\[
\theta(\bullet,\bar{y})^{\prime}(\bar{x};v) \, = \,
\displaystyle{
\lim_{\nu \to \infty}
} \, \displaystyle{
\frac{\theta(x^{\, \nu},\bar{y}) - \theta(\bar{x},\bar{y})}{\tau_{\nu}}
} \, = \, \theta^{\, \prime}((\bar{x},\bar{y});(v,0)) \, \geq \, 0,
\]
where the two equalities hold by the B-differentiability of $\theta$.
\end{proof}

Let
\begin{equation} \label{eq:B-stat tuple epi}
\left\{  \, \bar{x}; \, \{ \, \bar{t}_k \}_{k=1}^K; \, \{ \, \bar{s}_{\ell} \, \}_{\ell=1}^L \, \right\}
\end{equation}
be a B-stationary tuple of (\ref{eq:original POP epi form}).  By Lemma~\ref{lm:stationary of bivariate},
it follows that the tuple $\left\{ \, \{ \, \bar{t}_k \}_{k=1}^K; \, \{ \, \bar{s}_{\ell} \, \}_{\ell=1}^L \, \right\}$
is a B-stationary solution of the problem:
\[ \begin{array}{ll}
\displaystyle{
\operatornamewithlimits{\mbox{\bf minimize}}_{t; \, s}
} & \displaystyle{
\sum_{k=1}^K
} \, t_k + \lambda \,
\max\left( \, \displaystyle{
\sum_{\ell=1}^L
} \, s_{\ell} - b, \ 0 \, \right) \\ [0.2in]
\mbox{\bf subject to} &
t_k \, \geq \, \varphi_k(\bar{x}) \, \onebld_{( \, 0,\infty \, )}(g_k(\bar{x})), \epc \forall \, k \, \in \, [ K ] \\ [0.1in]
\mbox{and} & s_{\ell} \, \geq \, \phi_{\ell}(\bar{x}) \, \onebld_{( \, 0,\infty \, )}(h_{\ell}(\bar{x})), \epc \forall \, \ell \, \in \, [ L ].
\end{array}
\]
This is a trivial convex piecewise linear program in the $(t,s)$-variables with lower-bound constraints only.  It can
easily be seen that
$\bar{t}_k = \pi_k^{\, \varphi}(\bar{x})$ for all $k \in [ K ]$, but similar equalities are
not guaranteed for the $s$-variables.  In fact, there are two possibilities:
(i) $\displaystyle{
\sum_{\ell=1}^L
} \, \bar{s}_{\ell} \leq b$ and there exists $\bar{\ell} \in [ L ]$ such that
$\bar{s}_{\bar{\ell}} > \pi_{\bar{\ell}}^{\, \phi}(\bar{x})$; or
(ii) $\bar{s}_{\ell} = \pi_{\ell}^{\, \phi}(\bar{x})$ for all $\ell \in [ L]$.
The following analysis addresses both cases.  In case (i), it follows that
\[
\displaystyle{
\sum_{\ell=1}^L
} \, \pi_{\ell}^{\, \phi}(\bar{x}) \, < \, b,
\]
thus in particular $\bar{x}$ is feasible to (\ref{eq:original POP}).
We will return to complete this case in the main
Theorem~\ref{th:fixed-point stationary}.  For now, we consider the second case where both sets
of equalities in (\ref{eq:variable substitution}) hold at $\bar{x}$.
We next show that the vector $\bar{x}$ in the B-stationary tuple (\ref{eq:B-stat tuple epi})
must be a B-stationary solution of
the following problem defined with respect to a given $\bar{x} \in X$ and without the auxiliary variables $(t,s)$:
\[
\displaystyle{
\operatornamewithlimits{\mbox{\bf minimize}}_{x \in X}
} \ \wh{\Phi}_{\lambda}(x;\bar{x}) \, \triangleq \, c(x) +
\displaystyle{
\sum_{k \in {\cal K}_>(\bar{x})}
} \, \varphi_k(x) + \lambda \, \max\left( \, \displaystyle{
\sum_{\ell \in {\cal L}_>(\bar{x})}
} \, \phi_{\ell}(x) - b, \ 0 \, \right)
\]
\begin{equation} \label{eq:POP epi pulled-down}
\left. \begin{array}{ll}
\mbox{\bf subject to} & g_k(x) \, \leq \, 0 \epc \forall \, k \in {\cal K}_{\leq}(\bar{x}) \\ [0.1in]
& g_k(x) \, \geq \, 0 \epc \forall \, \ell \in {\cal K}_>(\bar{x}) \\ [0.1in]
& h_{\ell}(x) \, \leq \, 0 \epc \forall \, \ell \in {\cal L}_{\leq}(\bar{x}) \\ [0.1in]
\mbox{\bf and} & h_{\ell}(x) \, \geq \, 0 \epc \forall \, \ell \in {\cal L}_>(\bar{x})
\end{array} \right\} \epc \mbox{\begin{tabular}{l}
including $x \in X$, the set of these constraints \\ [5pt]
is denoted by $\wh{S}_{\rm ps}(\bar{x})$ and is equal to \\ [5pt]
$S_{\rm ps}(\bar{x})$ without the functional constraint
\end{tabular}}
\end{equation}
Let $dx$ be a tangent vector of $\wh{S}_{\rm ps}(\bar{x})$ at $\bar{x}$.  There exist sequences $\{ x^{\, \nu} \} \to \bar{x}$ and
$\{ \tau_{\nu} \} \downarrow 0$ such that
$x^{\nu} \in \wh{S}_{\rm ps}(\bar{x})$ for all $\nu$ and $dx = \displaystyle{
\lim_{\nu \to \infty}
} \, \displaystyle{
\frac{x^{\, \nu} - \bar{x}}{\tau_{\nu}}
}$.  Define
\[
dt_k \, \triangleq \, \left\{ \begin{array}{cl}
\varphi_k^{\, \prime}(\bar{x};dx) & \mbox{if $k \in {\cal K}_>(\bar{x})$} \\ [5pt]
0 & \mbox{otherwise,}
\end{array} \right. \epc \mbox{and} \epc
ds_{\ell} \triangleq \left\{ \begin{array}{cl}
\phi_{\ell}^{\, \prime}(\bar{x};dx) & \mbox{if $\ell \in {\cal L}_>(\bar{x})$} \\ [5pt]
0 & \mbox{otherwise.}
\end{array} \right.
\]
We claim that
\[
( \, dt_k,dx \, ) \, \in \, {\cal T}({\rm epi}(\pi_k^{\varphi};X);(\bar{t}_k,\bar{x})) \ \mbox{ and } \
( \, ds_{\ell},dx \, ) \, \in \, {\cal T}({\rm epi}(\pi_{\ell}^{\phi};X);(\bar{s}_{\ell},\bar{x})), \epc
\forall \, ( k,\ell ) \in [ K ] \times [ L ].
\]
This is obviously true for all $k \in {\cal K}_>(\bar{x}) \cup {\cal K}_<(\bar{x})$ and $\ell \in {\cal L}_>(\bar{x}) \cup {\cal L}_<(\bar{x})$,
by the first two subcases of Proposition~\ref{pr:epigraph product}(b).   For an index $k \in {\cal K}_=(\bar{x})$;
we have $g_k(x^{\, \nu} ) \leq 0$ by (\ref{eq:POP epi pulled-down}).  With $t_k^{\, \nu} \triangleq \, \pi_k^{\varphi}(x^{\, \nu}) = 0$,
recalling $\bar{t}_k = \pi_k^{\varphi}(\bar{x}) = 0$,
we obtain $dt_k = 0 = \displaystyle{
\lim_{\nu \to \infty}
} \, \displaystyle{
\frac{t_k^{\, \nu} - \bar{t}_k}{\tau_{\nu}}
}$.  Therefore, $( \, dt_k,dx \, )$ belongs to ${\cal T}({\rm epi}(\pi_k^{\varphi};X);(\bar{t}_k,\bar{x}))$ for all $k \in [ K ]$.
Similarly, we can prove $( \, ds_{\ell},dx \, ) \in {\cal T}({\rm epi}(\pi_{\ell}^{\phi};X);(\bar{s}_{\ell},\bar{x}))$
for all $\ell \in [ L ]$.  Furthermore, it is easy to see that
$\Phi_{\lambda}^{\, \prime}((\bar{x},\bar{t},\bar{s});(dx,dt,ds)) = \wh{\Phi}_{\lambda}(\bullet;\bar{x})^{\, \prime}(\bar{x},dx)$.  This
completes the proof that $\bar{x}$ is a B-stationary solution of (\ref{eq:POP epi pulled-down}).

\gap

So far, we have not imposed any constraint qualifications; for the last step in concluding that $\bar{x}$ is B-stationary solution
of (\ref{eq:barx stationary problem}), we need the last condition imposed in Theorem\ref{th:fixed-point stationary} below.
This condition is the one
stated in \cite[Theorem~9.2.1]{CuiPang2021} specialized to the problem (\ref{eq:POP epi pulled-down}); we refer to this reference
for a brief history of the condition in the theory of exact penalization.  For the problem (\ref{eq:original POP no Heaviside constraint})
this extra condition is not needed; therefore we obtain a constructive approach for obtaining a pseudo B-stationary solution of
this problem, which has no constraint involving the Heaviside functions,
under only the sign conditions on the $\varphi_k$ functions (and the blanket B-differentiability assumption of the problem).
 {For the sake of clarity in the last step of the proof of Theorem~\ref{th:fixed-point stationary},
we state the following simple lemma but omit its proof.}

\begin{lemma} \label{lm:dd of max composite} \rm
Let $f : {\cal O} \to \mathbb{R}$ be B-differentiable at $\bar{x}$.  Let $f_{\max}(x) \triangleq \max( f(x), 0)$.  Then
\[
f_{\max}^{\, \prime}(\bar{x};v) \, \leq \, \max\left( \, f^{\, \prime}(\bar{x};v), \, 0 \, \right), \epc \forall \, v \, \in \, \mathbb{R}^n.
\]
\end{lemma}


\begin{theorem} \label{th:fixed-point stationary} \rm
Under the blanket assumption of problem (\ref{eq:original POP}), assume that $\varphi_k$ and $\phi_{\ell}$
are nonnegative on $X \cap g_k^{-1}(0)$ and $X \cap h_{\ell}^{-1}(0)$, respectively, for all $( k,\ell ) \in [ K ] \times [ L ]$,
and that $c$ and each $\varphi_k$ are Lipschitz continuous on $X$ with Lipschitz constants $\mbox{Lip}_c$ and
$\mbox{Lip}_{\varphi}$, respectively.  If $\left\{ \bar{x}, \bar{t}, \bar{s} \right\}$
is a B-stationary tuplet of (\ref{eq:original POP epi form}) corresponding to a $\lambda$ satsfying
\begin{equation} \label{eq:condition on penalty}
\lambda \, > \, \mbox{Lip}_c + K \, \mbox{Lip}_{\varphi},
\end{equation}
then $\bar{x}$ is a pseudo B-stationary solution of the problem (\ref{eq:original POP}), provided that there exists a vector
$\bar{v} \in {\cal T}( \wh{S}_{\rm ps}(\bar{x});\bar{x})$ with unit length satisfying:
$\displaystyle{
\sum_{\ell \in {\cal L}_>(\bar{x})}
} \, \phi_{\ell}^{\, \prime}(\bar{x};\bar{v}) \leq -1$.
\end{theorem}

\begin{proof}  Continuing the above analysis of the case where (\ref{eq:variable substitution}) holds
at $\bar{x}$, we suppose that $\bar{x}$ fails the functional constraint; i.e.,
$\displaystyle{
\sum_{\ell=1}^L
} \, \pi_{\ell}^{\phi}(\bar{x}) > b$.  Then we have
\[ \begin{array}{lll}
0 & \leq & \wh{\Phi}_{\lambda}(\bullet;\bar{x})^{\, \prime}(\bar{x},\bar{v}) \\ [0.1in]
& = & c^{\, \prime}(\bar{x};\bar{v}) + \displaystyle{
\sum_{k \in {\cal K}_>(\bar{x})}
} \, \varphi_k^{\, \prime}(\bar{x};\bar{v}) + \lambda \, \displaystyle{
\sum_{\ell \in {\cal L}_>(\bar{x})}
} \, \phi_{\ell}^{\, \prime}(\bar{x};\bar{v}) \\ [0.2in]
& \leq & ( \, \mbox{Lip}_c + K \, \mbox{Lip}_{\varphi} \, ) - \lambda .
\end{array} \]
This contradicts the condition (\ref{eq:condition on penalty}).  Hence $\bar{x}$ is feasible to (\ref{eq:original POP}).
To show that $\bar{x}$ is B-stationary for (\ref{eq:barx stationary problem}), there are two cases to consider: $\displaystyle{
\sum_{\ell \in {\cal L}_>(\bar{x})}
} \, \phi_{\ell}(\bar{x}) = b$ or $\displaystyle{
\sum_{\ell \in {\cal L}_>(\bar{x})}
} \, \phi_{\ell}(\bar{x}) < b$.  Consider the former case first.
Let $v \in {\cal T}(S_{\rm ps}(\bar{x});\bar{x})$.
Then $v \in {\cal T}(\wh{S}_{\rm ps}(\bar{x});\bar{x})$ and $\displaystyle{
\sum_{\ell \in {\cal L}_>(\bar{x})}
} \, \phi_{\ell}^{\, \prime}(\bar{x};v) \leq 0$.
We have
\begin{equation} \label{eq:stationarity in proof}
0 \, \leq \, \wh{\Phi}_{\lambda}(\bullet;\bar{x})^{\, \prime}(\bar{x},v) \, = \,
c^{\, \prime}(\bar{x};v) + \displaystyle{
\sum_{k \in {\cal K}_>(\bar{x})}
} \, \varphi_k^{\, \prime}(\bar{x};v),
\end{equation}
establishing the B-stationarity of $\bar{x}$ for the problem (\ref{eq:barx stationary problem}).  In the latter case where
$\displaystyle{
\sum_{\ell \in {\cal L}_>(\bar{x})}
} \, \phi_{\ell}(\bar{x}) < b$, we have ${\cal T}(S_{\rm ps}(\bar{x});\bar{x}) = {\cal T}(\wh{S}_{\rm ps}(\bar{x});\bar{x})$
and the same expression (\ref{eq:stationarity in proof}) also holds.

\gap

What remains to be proved is the case where the B-stationary tuple (\ref{eq:B-stat tuple epi}) of (\ref{eq:original POP epi form}) is
such that
\[ \displaystyle{
\sum_{\ell \in {\cal L}_>(\bar{x})}
} \, \phi_{\ell}(\bar{x}) \, = \,
\displaystyle{
\sum_{\ell=1}^L
} \, \pi_{\ell}^{\, \phi}(\bar{x}) \, < \, \displaystyle{
\sum_{\ell=1}^L
} \, \bar{s}_{\ell} \, \leq \, b
\]
and $\bar{s}_{\bar{\ell}} > \pi_{\bar{\ell}}^{\phi}(\bar{x})$ for at least one $\bar{\ell}$.
Let $dx \in {\cal T}(S_{\rm ps}(\bar{x});\bar{x})$ be arbitrary.  There exist sequences $\{ x^{\, \nu} \} \to \bar{x}$ and
$\{ \tau_{\nu} \} \downarrow 0$ such that
$x^{\nu} \in S_{\rm ps}(\bar{x}) \subseteq \wh{S}_{\rm ps}(\bar{x})$ for all $\nu$ and $dx = \displaystyle{
\lim_{\nu \to \infty}
} \, \displaystyle{
\frac{x^{\, \nu} - \bar{x}}{\tau_{\nu}}
}$.  As we have already shown,
$( dt_k,dx ) \in {\cal T}({\rm epi}(\pi_k^{\varphi};X);(\bar{t}_k,\bar{x}))$
for all $k \in [ K ]$, where $dt_k \triangleq \left\{ \begin{array}{cl}
\varphi_k^{\, \prime}(\bar{x};dx) & \mbox{if $k \in {\cal K}_>(\bar{x})$} \\ [5pt]
0 & \mbox{otherwise.}
\end{array} \right.$
For an arbitrary scalar $M > 0$, define a vector $ds^M$ as follows:
\[
ds_{\ell}^M \, \triangleq \left\{ \begin{array}{cl}
\phi_{\ell}^{\, \prime}(\bar{x};dx) & \mbox{if $\ell \in {\cal L}_>(\bar{x})$ and $\bar{s}_{\ell} = \pi_{\ell}^{\phi}(\bar{x})$}, \epc
\mbox{index set denoted ${\cal L}_>^=(\bar{x})$} \\ [0.1in]
-M & \mbox{if $\bar{s}_{\ell} > \pi_{\ell}^{\phi}(\bar{x})$} \\ [0.1in]
0 & \mbox{otherwise.}
\end{array} \right.
\]
As before, we have $( ds_{\ell}^M,dx ) \in {\cal T}({\rm epi}(\pi_{\ell}^{\phi};X);(\bar{s}_{\ell},\bar{x}))$
for all $\ell \in [ L ]$ except possibly when $\bar{s}_{\ell} > \pi_{\ell}^{\phi}(\bar{x})$.
We show that the latter exception can be removed.  Indeed, with $M$ fixed but arbitrary, for an index $\ell$ of the
latter kind, there are 2 cases to consider:

\gap

(i) $h_{\ell}(\bar{x}) \neq 0$: then
$( -M,dx ) \in {\cal T}({\rm epi}(\pi_{\ell}^{\phi};X);(\bar{s}_{\ell},\bar{x}))$ by part (c) of Proposition~\ref{pr:epigraph product}.
%
%

\gap

(ii) $h_{\ell}(\bar{x}) = 0$: then $\bar{s}_{\ell} > \pi_{\ell}^{\, \phi}(\bar{x}) = 0 \geq h_{\ell}(x^{\, \nu})$,
where the last inequality holds because $x^{\nu} \in S_{\rm ps}(\bar{x})$.  Thus,
$\bar{s}_{\ell} - \tau_{\nu} M > \pi_{\ell}^{\, \phi}(x^{\, \nu}) = 0$
for all $\nu$ sufficiently large and $( -M,dx ) \in {\cal T}({\rm epi}(\pi_{\ell}^{\phi};X);(\bar{s}_{\ell},\bar{x}))$ follows.

\gap

We have thus completed the proof that $( ds_{\ell}^M,dx ) \in {\cal T}({\rm epi}(\pi_{\ell}^{\phi};X);(\bar{s}_{\ell},\bar{x}))$
for all $\ell \in [ L ]$.  We have
\[
\displaystyle{
\sum_{\ell=1}^L
} \, ds_{\ell}^M \, = \, \displaystyle{
\sum_{\ell \in {\cal L}_>^=(\bar{x})}
} \, \phi_{\ell}^{\, \prime}(\bar{x};dx) - M \, \underbrace{| \, \{ \, \ell \, : \, \bar{s}_{\ell} > \pi_{\ell}^{\phi}(\bar{x}) \, \} \, |}_{\mbox{$\geq 1$}}.
\]
By choosing $M$ sufficiently large, we have $\displaystyle{
\sum_{\ell=1}^L
} \, ds_{\ell}^M \leq 0$.  By Lemma \ref{lm:dd of max composite}, it follows that
\[ \begin{array}{lll}
0 & \leq & \Phi_{\lambda}( (\bar{x},\bar{t},\bar{s});(dx,dt,ds^M) ) \\ [0.1in]
& \leq & c^{\, \prime}(\bar{x};dx) + \displaystyle{
\sum_{k \in {\cal K}_>(\bar{x})}
} \, \varphi_k^{\, \prime}(\bar{x};dx) + \lambda \, \max\left( \,  {\displaystyle{
\sum_{\ell}^L
} \, ds_{\ell}^M}, \, 0 \, \right) \\ [0.2in]
&  {=} &  c^{\, \prime}(\bar{x};dx) + \displaystyle{
\sum_{k \in {\cal K}_>(\bar{x})}
} \, \varphi_k^{\, \prime}(\bar{x};dx),
\end{array} \]
showing that $\bar{x}$ is a B-stationary solution of (\ref{eq:barx stationary problem}), as desired.
\end{proof}

\section{A Digression: Approximations of the (open) Heaviside Function} \label{sec:approximations Heaviside}

The lifted formulation (\ref{eq:original POP epi form}), or its computational workhorse (\ref{eq:original POP exposed epi form}),
 {requires the auxiliary variables $t$ and $s$ but has the advantage of leading}
directly to a pseudo B-stationary point of the original problem (\ref{eq:original POP}) under some
mild conditions as stated in Theorem~\ref{th:fixed-point stationary}.  Of independent interest,
the next approach is based on approximations of the open Heaviside function and leads to approximated problems in the $x$-variable only.
In this section, we digress from the discussion of the problem (\ref{eq:original POP}) and present a focused discussion on
two approaches to construct such approximations and show how they are related to each other.  One construction is based on truncation \cite{CuiLiuPang21}
and the other is based on the classical work of \cite{ENWets95} on mollifiers, or smoothing \cite{Chen12}.

\begin{definition} \label{def:approximation} \rm
We say that the bivariate functions $\theta : \mathbb{R} \times \mathbb{R}_{++} \to [ \, 0,1 \, ]$
{\sl p-approximate} (``p'' for pointwise) the (open) Heaviside function
$\onebld_{( \, 0,\infty \, )}$ 
if there exist endpoint functions $\underline{\theta}$ and $\overline{\theta} : \mathbb{R}_{++} \to \mathbb{R}_+$ satisfying

\gap

(A0) $\displaystyle{
\lim_{\delta \downarrow 0}
} \, \underline{\theta}(\delta) \, = \, \displaystyle{
\lim_{\delta \downarrow 0}
} \, \overline{\theta}(\delta) \, = \, 0$;

\gap

(A1) $\displaystyle{
\lim_{\delta \downarrow 0}
} \, \theta(t,\delta) = \onebld_{( \, 0,\infty \, )}(t)$ for all $t \in \mathbb{R}$; [this condition allows us to define
$\theta(t,0) \triangleq \onebld_{( \, 0,\infty \, )}(t)$, thereby extending the domain of definition of $\theta$ to
$\mathbb{R} \times \mathbb{R}_+$]; and

\gap

(A2) $\theta(t,\delta) = \left\{ \begin{array}{ll}
0 & \forall \ t \, \leq \, -\underline{\theta}(\delta) \\ [5pt]
1 & \forall \ t \, \geq \ \overline{\theta}(\delta)
\end{array} \right\}$ for all $\delta > 0$.

\gap

For subsequent purposes, we further stipulate that for all $\delta > 0$,

\gap

(A3) the function $\theta(\bullet,\delta)$ is B-differentiable on $[ \, -\underline{\theta}(\delta), \, \overline{\theta}(\delta) \, ]$
with the directional derivatives satisfying
\begin{equation} \label{eq:A3 for p-approx}
\theta(\bullet,\delta)^{\, \prime}(t;1) \, \geq \, 0 \, \geq \,
\theta(\bullet,\delta)^{\, \prime}(t;-1), \epc \forall \, t \, \in \, [ \, -\underline{\theta}(\delta), \, \overline{\theta}(\delta) \, ].
\end{equation}
\end{definition}

By (A2), condition (A3) implies that $\theta(\bullet,\delta)^{\, \prime}(t;1) \, \geq \, 0 \, \geq \,
\theta(\bullet,\delta)^{\, \prime}(t;-1)$ for all $t \in \mathbb{R}$.
Clearly, if $\theta$ p-approximates the Heaviside function $\onebld_{( \, 0,\infty \, )}$, then so does the function
$\psi(t,\delta) \triangleq \theta\left( \displaystyle{
\frac{t}{m(\delta)}
}, \ \delta \right)$ with the variable $t$ being scaled by the positive univariate
function $m : \mathbb{R}_{++} \to \mathbb{R}_{++}$ satisfying $\displaystyle{
\limsup_{\delta \downarrow 0}
} \, m(\delta) < \infty$.
In the next two subsections, we present two ways to obtain p-approximiations of the Heaviside function satisfying the B-differentiability requirement.

\subsection{Truncation derived approximations}

Let $\wh{\theta} : \mathbb{R} \times \mathbb{R}_{++} \to \mathbb{R}$
be such that there exist end-point functions $\underline{\theta}$ and $\overline{\theta} : \mathbb{R}_{++} \to \mathbb{R}_+$ satisfying (A0)
and

\gap

(T1) $\wh{\theta}(t,\delta) \, \left\{ \begin{array}{ll}
\leq \, 0 & \forall \ t \, \leq \, -\underline{\theta}(\delta) \\ [5pt]
\, \geq \, 1 & \forall \ t \, \geq \ \overline{\theta}(\delta)
\end{array} \right\}$ for all $\delta > 0$;

\gap

(T2) $\displaystyle{
\lim_{\delta \downarrow 0}
} \, \wh{\theta}(0,\delta) = 0$;

\gap

(T3) $\wh{\theta}(\bullet,\delta)$ is B-differentiable on an open interval containing
$[ \, -\underline{\theta}(\delta),\overline{\theta}(\delta) \, ]$
with
\[
\wh{\theta}(\bullet,\delta)^{\, \prime}(t;1) \, \geq \, 0 \, \geq \wh{\theta}(\bullet,\delta)^{\, \prime}(t;-1), \epc \forall \,
t \, \in \, [ \, -\underline{\theta}(\delta),\overline{\theta}(\delta) \, ].
\]
To obtain the p-approximation function $\theta$ from $\wh{\theta}$, let
\[
T_{[ \, 0,1 \, ]}(t) \, \triangleq \, \min\left\{ \, \max( t,0 ), \, 1 \, \right\}
\, = \, \max\left\{ \, \min( t,1 ), \, 0 \, \right\} 
\, = \, \max( t,0 ) - \max( t-1, 0), \epc t \, \in \, \mathbb{R}
\]
be the truncation operator to the range $[ \, 0,1 \, ]$ and define the composite function:
\begin{equation} \label{eq:approx by truncation}
\theta_{\rm tr}(t,\delta) \, \triangleq \, T_{[ \, 0,1 \, ]}\left( \, \wh{\theta}(t,\delta) \, \right), \epc
( t,\delta ) \, \in \, \mathbb{R} \times \mathbb{R}_{++}.
\end{equation}
We formally state that the truncated function $\theta$ is a p-approximation of the Heaviside function
in the first part of the next proposition.  The truncation function was used in \cite{CuiLiuPang21} as a unification scheme for many
approximations of chance constraints in stochastic programs; see the cited reference for many prior works on the latter subject.
Part (b) of the result provides easy sufficient conditions for (T1) and (T3) to hold; this part is the bridge between the truncation
approach discussed herein and the nonifier approach to be discussed in the next subsection.  Part (c) presents a broad family of approximating
functions by truncation
that generalize the perspective functions studied extensively in convex
analysis \cite{BauschkeCombettes11,Combettes18,HULemarechal93,rockafellar1970convex}
and employed extensively for integer programs in recent years; see e.g.\ \cite{GunlukLinderoth10}.  The approximating functions
in the last part of the proposition yield the family of folded concave approximations of the $\ell_0$ function \cite[Subsection~3.1.4]{CuiPang2021};
these include the smoothly clipped absolute deviation {\sc scad} function \cite{FanLi01}; the minimax concave
penalty {\sc mcp} function \cite{Zhang10,DongChenLinderoff15}; the capped $\ell_1$-function \cite[Section~5]{LeThiTaoVo15},
and others.  As shown in \cite{AhnPangXin17}, all these functions are of the difference-of-convex kind and
not differentiable at the origin as they are approximations of the $\ell_0$-function that is discontinuous there.

\begin{proposition} \label{pr:approx by truncation} \rm
The following statements hold:

\gap

(a) If $\wh{\theta}$ satisfies (T1), (T2), and (T3), then
its truncation function $\theta_{\rm tr}$ p-approximates the
Heaviside function $\onebld_{( \, 0,\infty \, )}$ and satisfies (A3).

\gap

(b) If $\wh{\theta}(\bullet,\delta)$ is nondecreasing and satisfies:

\gap

(T1$^{\, \prime}$) $\wh{\theta}(-\underline{\theta}(\delta),\delta) = 0$ and $\wh{\theta}(\overline{\theta}(\delta),\delta) = 1$ for all $\delta > 0$,

\gap

then (T1) holds; if additionally $\wh{\theta}(\bullet,\delta)$ is locally Lipschitz, then (T3) holds.

\gap

(c) If $\psi : \mathbb{R} \to \mathbb{R}$ is a nondecreasing B-differentiable function satisfying
$\psi(0) = 0$ and $\psi(1) = 1$, then
\[
\theta(t,\delta) \, \triangleq \, T_{[ \, 0,1 \, ]}\left( \psi\left( q(\delta) + \displaystyle{
\frac{t}{m(\delta)}
} \, \right) \, \right)
\]
with $q : \mathbb{R}_{++} \to [ \, 0,1 \, ]$ and
$m : \mathbb{R}_{++} \to \mathbb{R}_{++}$ satisfying
\begin{equation} \label{eq:limits on q and m}
\displaystyle{
\lim_{\delta \downarrow 0}
} \, q(\delta) \, = \, 0 \, = \, \displaystyle{
\lim_{\delta \downarrow 0}
} \, m(\delta) 
\end{equation}
p-approximates $\onebld_{( \, 0,\infty \, )}$ and satisfies (A3).

\gap

(d) If $\theta_1(t,\delta)$ and $\theta_2(t,\delta)$ both p-approximate the open Heaviside function $\onebld_{( \, 0,\infty \, )}$,
then their sum $\theta_1(t,\delta) + \theta_2(-t,\delta)$ p-approximates the $\ell_0$-function $| \, t \, |_0$.
\end{proposition}

\begin{proof}  For (a),
we need to show that (T1) and (T3) imply (A1) and (A2) and that (T3) implies (A3).  Clearly, (T2) implies (A1) for $t = 0$.
For a $t_* > 0$, we have
$t_* \geq \overline{\theta}(\delta)$ for all $\delta > 0$ sufficiently small.  Thus,
$\theta_{\rm tr}(t_*,\delta) = 1 = \onebld_{( \, 0,\infty \, )}(t_*)$ by (T1) and the definition of truncation, which is (A1) at $t_*$.
Similarly, (A1) also holds for $t_* < 0$.  Thus the pointwise convergence condition (A1) holds for all $t \in \mathbb{R}$.
Clearly (A2) holds by truncation.  For (A3), we have
\[
\theta_{\rm tr}(\bullet,\delta)^{\, \prime}(t;\pm 1) \, = \, T_{[ \, 0,1 \, ]}^{\, \prime}( \wh{\theta}(t,\delta);
\wh{\theta}(\bullet,\delta)^{\, \prime}(t;\pm 1) ).
\]
Since $T_{[ \, 0,1 \, ]}^{\, \prime}(t;1) \geq 0 \geq T_{[ \, 0,1 \, ]}^{\, \prime}(t;-1)$ for all $t \in \mathbb{R}$,
we readily obtain (A3) from (T3).  If $\wh{\theta}$ is nondecreasing, then (T1$^{\, \prime}$)
clearly implies (T1); moreover, this function is directionally differentiable and by the definition (\ref{eq:dd of univariate}),
we have
$\wh{\theta}(\bullet,\delta)^{\, \prime}(t;1) \geq 0 \geq \wh{\theta}(\bullet,\delta)^{\, \prime}(t;-1)$
for all $t$; statement (b) thus holds.
For statement (c), define $\underline{\theta}(\delta) \triangleq -m(\delta) \, q(\delta)$ and
$\overline{\theta}(\delta) \triangleq m(\delta) \, ( 1 - q(\delta) )$.  By the limit on the function $m$ in (\ref{eq:limits on q and m}),
we obtain $\displaystyle{
\lim_{\delta \downarrow 0}
} \, \overline{\theta}(\delta) = 0 = \displaystyle{
\lim_{\delta \downarrow 0}
} \, \underline{\theta}(\delta)$.  Finally, with $\wh{\theta}(t,\delta) \triangleq \psi\left( q(\delta) + \displaystyle{
\frac{t}{m(\delta)}
} \, \right)$, condition (T2) holds also by the limits on $q(\delta)$ in  (\ref{eq:limits on q and m}).  The last statement (d)
is obvious by (\ref{eq:ell0 fnc}).
\end{proof}

To illustrate the function $\theta$ in part (c) of the above proposition and the role of condition (T2)
consider a common approximation of the Heaviside function
$\onebld_{( \, 0,\infty \, )}$ by the truncated hinge loss functions \cite{WuLiu07,QiCuiLiuPang19}: $T_{\rm h}(t,\delta) \triangleq \displaystyle{
\frac{1}{2 \delta}
} \, \left[ \, \max( t + \delta, 0 ) - \max( t - \delta, 0 ) \, \right]$ for $\delta > 0$.  Notice that $T_{\rm h}(0,\delta) = \thalf$ for all $\delta$;
thus this approximation function fails condition (T2) and as a result does not ``recover'' the Heaviside function as $\delta \downarrow 0$.
When this function $T_{\rm h}(\bullet,\delta)$ is
employed in sampled discretization of a probability function, as in the cited references:
\[
\mathbb{P}_{\tilde{z}}( f(x,\tilde{z}) > 0 ) \, = \, \mathbb{E}_{\tilde{z}}[ \onebld_{( \, 0,\infty \, )}( f(x,\tilde{z}) ) ] \, \approx \,
\displaystyle{
\frac{1}{N}
} \, \displaystyle{
\sum_{s=1}^N
} \, T_{\rm h}( f(x,z^s),\delta ),
\]
where $\{ z^{\, s} \}_{s=1}^N$ is a sample batch of size $N$ of the random variable $\tilde{z}$,
the gap between the indicator function and its approximation at the origin is less important because of the common assumption that
$\mathbb{P}_{\tilde{z}}( f(x,\tilde{z}) = 0 ) = 0$.  Nevertheless, in a deterministic context, such a gap can be significant because the set $f^{-1}(0)$
is typically of most interest for the composite function $\onebld_{( \, 0,\infty \, )}(f(x))$.  Part of the reason for the gap is due to the
symmetry of the function $T_{\rm h}(\bullet,\delta)$ with respect to $t = 0$.  To recover the Heaviside function exactly as $\delta \downarrow 0$, we may
consider the following modified hinge loss function:
\begin{equation} \label{eq:nonifier for modified hinge}
\wt{T}_{\rm h}(t,\delta) \, \triangleq \, \min\left\{ \, \max\left( \, \displaystyle{
\frac{t}{\delta + \sqrt{\delta}}
} + \displaystyle{
\frac{\sqrt{\delta}}{1 + \sqrt{\delta}}
} \, , \ 0 \ \right), \ 1 \, \right\}, \epc \mbox{for $\onebld_{( \, 0, \infty \, )}(t)$},
\end{equation}
which is derived from $\psi$ being the identity function, $q(\delta) = \displaystyle{
\frac{\sqrt{\delta}}{1 + \sqrt{\delta}}
}$ and $m(\delta) = \delta + \sqrt{\delta}$, 
both satisfying the limits (\ref{eq:limits on q and m}).

\subsection{Nonifiers induced approximations}

We present the other approach for deriving approximations for the Heaviside function based on
``averaged functions'' as defined in \cite[Definition~3.1]{ENWets95}.

\begin{definition} \label{df:averaged functions} \rm
Given a locally integrable function $f : \mathbb{R}^n \to \mathbb{R}$ and a family of bounded mollifiers $\{ \psi(\bullet,\delta) : \mathbb{R}^n \to \mathbb{R}_+;
\, \delta \in \mathbb{R}_+ \}$ that satisfy
\[ \displaystyle{
\int_{\mathbb{R}^n}
} \, \psi(z,\delta) \, dz \, = \, 1, \epc \mbox{supp } \psi(\bullet,\delta) \, \triangleq \, \{ \, z \, \in \, \mathbb{R}^n \, \mid \,
\psi(z,\delta) \, > \, 0 \, \}
\subseteq \ \rho_{\delta} \, \ball, \mbox{ for some } \{ \rho_{\delta} \} \downarrow 0 \mbox{ as } \delta \downarrow 0,
\]
where $\ball$ is a unit Euclidean ball in $\mathbb{R}^n$, the associated family $\{ f_{\psi}(\bullet,\delta) \}_{\delta \geq 0}$ of {\sl averaged functions}
is given by
\[
f_{\psi}(x,\delta) \, \triangleq \, \displaystyle{
\int_{\mathbb{R}^n}
} \, f(x - z) \, \psi(z,\delta) \, dz \, = \, \displaystyle{
\int_{\mathbb{R}^n}
} \, f(z) \, \psi(x - z,\delta) \, dz.
\]
\end{definition}

As noted in \cite{ENWets95}, $\psi(\bullet,\delta)$ is a probability density function whose support tends to zero as $\delta \downarrow 0$;
moreover, $\psi(\bullet,\delta)$ does not need to be continuous.
Focus of the references \cite{ENWets95,Chen12} has been on mollifiers $\psi(\bullet,\delta)$ that lead to smooth (i.e., continuously differentiable) averaged functions $f_{\psi}(\bullet,\delta)$.  Two sets of conditions ensure the latter property: (i) $f$ is continuous with certain special choices of the
family $\{ \psi(\bullet,\delta) \}$
(see \cite[Proposition~3.11]{ENWets95}), and (ii) the mollifiers $\psi(\bullet,\delta)$ are continuously differentiable (see \cite[Proposition~3.9]{ENWets95}).
Since we are most interested in nonsmooth averaged functions, we coin the term ``nonifier'' for $\psi$ with the intention that
the induced averaged functions $f_{\psi}(\bullet,\delta)$ are not necessarily differentiable.
In what follows, we derive these averaged functions for the Heaviside function and show how they are related to those obtained from
the previous truncation approach.

\gap

To begin, we take a bivariate function $\psi : \mathbb{R} \times \mathbb{R}_{++} \to \mathbb{R}_+$ such that $\psi(\bullet,\delta)$ is integrable
on $\mathbb{R}$ for every $\delta > 0$ with $\displaystyle{
\int_{-\infty}^{\, \infty}
} \, \psi(t,\delta) \, dt = 1$ and there exist
end-point functions $\underline{\psi}$ and
$\overline{\psi} : \mathbb{R}_{++} \to \mathbb{R}_+$ satisfying condition (A0), i.e., $\displaystyle{
\lim_{\delta \downarrow 0}
} \, \underline{\psi}(\delta) = 0 = \displaystyle{
\lim_{\delta \downarrow 0}
} \, \overline{\psi}(\delta)$, such that
$\mbox{supp}(\psi(\bullet,\delta)) \subseteq [ \, -\underline{\psi}(\delta), \, \overline{\psi}(\delta) \, ]$.
It then follows that
\[ \begin{array}{lll}
\theta_{\, \psi}(t,\delta) & = & \displaystyle{
\int_0^{\, \infty}
} \, \psi(t - s,\delta) \, ds \, = \, \displaystyle{
\int_{-\infty}^{\, t}
} \, \psi(s,\delta) \, ds, \epc \forall \, ( \, t,\delta \, ) \, \in \, \mathbb{R} \times \mathbb{R}_{++} \\ [0.2in]
& = & \displaystyle{
\int_{-\underline{\psi}(\delta)}^{\min(t,\overline{\psi}(\delta))}
} \, \psi(s,\delta) \, ds \, = \, \left\{ \begin{array}{ll}
\epc 1 & \mbox{if $t \geq \overline{\psi}(\delta)$} \\ [0.1in]
\displaystyle{
\int_{-\underline{\psi}(\delta)}^{\, t}
} \, \psi(s,\delta) \, ds & \mbox{if $-\underline{\psi}(\delta) \leq t \leq \overline{\psi}(\delta)$} \\ [0.2in]
\epc 0 & \mbox{if $t \leq -\underline{\psi}(\delta)$}.
\end{array} \right.
\end{array} \]
Defining the cumulative distribution function:
\begin{equation} \label{eq:nonifier before truncation}
\wh{\theta}_{\, \psi}(t,\delta) \, \triangleq \, \displaystyle{
\int_{-\underline{\psi}(\delta)}^{\, t}
} \, \psi(s,\delta) \, ds, \epc \forall \, ( t,\delta ) \, \in \, \mathbb{R} \times \mathbb{R}_{++},
\end{equation}
we see that $\wh{\theta}_{\psi}(\bullet,\delta)$ is equal to its own truncation; i.e.,
$\wh{\theta}_{\, \psi}(t,\delta) = T_{[ \, 0,1 \, ]}(\wh{\theta}_{\, \psi}(t,\delta)) = \theta_{\, \psi}(t,\delta)$;
moreover, $\wh{\theta}_{\psi}(\bullet,\delta)$ is nondecreasing (because $\psi(\bullet,\delta)$ is nonnegative)
and satisfies condition (T1$^{\prime}$).  The following
lemma pertains to conditions (T2) and (T3) for the function $\wh{\theta}_{\, \psi}(\bullet,\delta)$;
no continuity of $\psi(\bullet,\delta)$ is needed.

\begin{lemma} \label{lm:properties nonifier open} \rm
The following two statements hold:

\gap

$\bullet $ If the one-sided limits $\psi(t\pm;\delta) \triangleq \displaystyle{
\lim_{\tau \downarrow 0}
} \, \psi(t \pm \tau,\delta)$ exist, then $\wh{\theta}_{\, \psi}(\bullet,\delta)$ is B-differentiable at $t$ with
\begin{equation} \label{eq:dd of integral}
\wh{\theta}_{\, \psi}(\bullet,\delta)^{\, \prime}(t;\pm 1) \, = \, \pm \, \psi(t\pm;\delta), \epc \mbox{respectively};
\end{equation}
thus condition (T3) (and (A3)) hold for $\wh{\theta}_{\, \psi}(\bullet,\delta)$; i.e.,
$\wh{\theta}_{\, \psi}(\bullet,\delta)^{\, \prime}(t;1) \geq 0 \geq
\wh{\theta}_{\, \psi}(\bullet,\delta)^{\, \prime}(t;-1)$.

\gap

$\bullet $ If $\displaystyle{
\limsup_{\delta \downarrow 0}
} \, \left\{ \, \left[ \, \underline{\psi}(\delta) + \overline{\psi}(\delta) \, \right] \, \displaystyle{
\sup_{t \in \mathbb{R}}
} \ \psi(t,\delta) \, \right\} < \infty$ and $\displaystyle{
\lim_{\delta \downarrow 0}
} \, \displaystyle{
\frac{\underline{\psi}(\delta)}{\overline{\psi}(\delta)}
} \, = \, 0$, then $\displaystyle{
\lim_{\delta \downarrow 0}
} \, \wh{\theta}_{\, \psi}(0,\delta) \, = \, 0$.  Thus condition (T2) holds for $\wh{\theta}_{\, \psi}(\bullet,\delta)$.


\end{lemma}

\begin{proof}  To show the locally Lipschtiz continuity of $\wh{\theta}_{\, \psi}(\bullet,\delta)$ near $t$, let $\varepsilon$ and $\bar{\tau}$ be
positive scalars such that
\[
\tau \, \in \, [ \, 0,\bar{\tau} \, ] \ \Rightarrow \ | \, \psi(t \pm \tau,\delta) - \psi(t\pm,\delta) \, | \, \leq \, \varepsilon.
\]
Let $t_1$ and $t_2$ be two scalars in $[ \, -\bar{\tau},\bar{\tau} \, ]$.  Consider first the case where both $t_1$ and $t_2$ are on the same
side of $t$.  We can write $t_1 = t + \tau_1$ and $t_2 = t + \tau_2$, with $\bar{\tau} \geq \tau_1 \geq \tau_2 \geq 0$.
\[
\wh{\theta}_{\, \psi}(t_1,\delta) - \wh{\theta}_{\, \psi}(t_2,\delta) \, = \, \displaystyle{
\int_{t_2}^{\, t_1}
} \, \psi(s,\delta) \, ds \, = \,  \displaystyle{
\int_{t_2}^{\, t_1}
} \, \left[ \, \psi(s,\delta) - \psi(t+,\delta) \, \right] \, ds + \psi(t+,\delta) \, ( \, t_1 - t_2 \, ).
\]
Hence,
\[
| \, \wh{\theta}_{\, \psi}(t_1,\delta) - \wh{\theta}_{\, \psi}(t_2,\delta) \, | \, \leq \, ( \, \varepsilon +  \psi(t+,\delta) \, ) \, ( \, t_1 - t_2 \, ).
\]
If $t_1$ an $t_2$ are on opposite sides of $t$, then can write $t_1 = t + \tau_1$ and $t_2 = t - \tau_2$ with $\tau_1$ and $\tau_2$ both
in the interval $[ \, 0,\bar{\tau} \, ]$.  It follows that
\[ \begin{array}{lll}
\wh{\theta}_{\, \psi}(t_1,\delta) - \wh{\theta}_{\, \psi}(t_2,\delta) & = &
 \displaystyle{
\int_{t_2}^{\, t}
} \, \psi(s,\delta) \, ds +  \displaystyle{
\int_t^{\, t_1}
} \, \psi(s,\delta) \, ds \\ [0.2in]
& = & \displaystyle{
\int_{t_2}^{\, t}
} \, \left[ \, \psi(s,\delta) - \psi(t-,\delta) \, \right] \, ds + \psi(t-,\delta) \, ( \, t - t_2 \, ) \ + \\ [0.2in]
& & \displaystyle{
\int_t^{\, t_1}
} \, \left[ \, \psi(s,\delta) - \psi(t+,\delta) \, \right] \, ds + \psi(t+,\delta) \, ( \, t_1 - t \, ) .
\end{array} \]
Hence,
\[ \begin{array}{lll}
| \, \wh{\theta}_{\, \psi}(t_1,\delta) - \wh{\theta}_{\, \psi}(t_2,\delta) \, | & \leq & ( \, \varepsilon +  \psi(t-,\delta) \, ) \, ( \, t - t_2 \, ) +
( \, \varepsilon +  \psi(t+,\delta) \, ) \, ( \, t_1 - t \, ) \\ [0.1in]
& \leq & \left[ \, \varepsilon +  \max\left( \, \psi(t-,\delta),\psi(t+,\delta) \, \right) \, \right] \,
[ \, ( \, t - t_2 \, ) + ( \, t_1 - t \, ) \, ] \\ [0.1in]
& = & \left[ \, \varepsilon +  \max\left( \, \psi(t-,\delta),\psi(t+,\delta) \, \right) \, \right] \, ( \, t_1 - t_2 \, ).
\end{array}
\]
The locally Lipschitz continuity of $\theta_{\psi}(\bullet,\delta)$ near $t$ follows.  Since $\theta_{\psi}(\bullet,\delta)$ is nondecreasing,
the B-differentiability of $\theta_{\psi}(\bullet,\delta)$ at $t$ follows.
For the two limits in (\ref{eq:dd of integral}), we prove only for the negative direction; i.e.,
$\wh{\theta}_{\, \psi}(\bullet,\delta)^{\, \prime}(t;-1) = -\psi(t-;\delta)$,
as the proof for the plus direction is similar (and a little more straightforward).  We have
\[ \begin{array}{l}
\wh{\theta}_{\, \psi}(\bullet,\delta)^{\, \prime}(t;-1) + \psi(t-,\delta) \, = \, \displaystyle{
\lim_{\tau \downarrow 0}
} \, \displaystyle{
\frac{\wh{\theta}_{\, \psi}(t - \tau,\delta) - \wh{\theta}_{\, \psi}(t,\delta) + \tau \, \psi(t-,\delta)}{\tau}
} \\ [0.2in]
\epc = \, \displaystyle{
\lim_{\tau \downarrow 0}
} \, \displaystyle{
\frac{\displaystyle{
\int_{t}^{\, t - \tau}
} \, \left[ \, \psi(s,\delta) - \psi(t-,\delta) \, \right] \, ds}{\tau}
} \, \leq \, \displaystyle{
\lim_{\tau \downarrow 0}
} \, \displaystyle{
\sup_{s \in [ \, t - \tau,t \, ]}
} \, | \, \psi(s,\delta) - \psi(t-,\delta) \, | \, = \, 0.
\end{array} \]
For the second statement, we have
\[ \begin{array}{lll}
\wh{\theta}_{\, \psi}(0,\delta) & = & \displaystyle{
\int_{-\underline{\psi}(\delta)}^{\, 0}
} \, \psi(s,\delta) \, ds \, \leq \, \underline{\psi}(\delta) \, \displaystyle{
\sup_{s \in \mathbb{R}}
} \, \psi(s,\delta) \, \leq \, \displaystyle{
\frac{\underline{\psi}(\delta)}{\underline{\psi}(\delta) + \overline{\psi}(\delta)}
} \, \left[ \, \underline{\psi}(\delta) + \overline{\psi}(\delta) \, \right] \, \displaystyle{
\sup_{t \in \mathbb{R}}
} \ \psi(t,\delta) \\ [0.2in]
& = & \displaystyle{
\frac{\underline{\psi}(\delta)/\overline{\psi}(\delta)}{1 + \underline{\psi}(\delta)/\overline{\psi}(\delta)}
} \, \left[ \, \underline{\psi}(\delta) + \overline{\psi}(\delta) \, \right] \, \displaystyle{
\sup_{t \in \mathbb{R}}
} \ \psi(t,\delta).
\end{array}
\]
Thus $\displaystyle{
\lim_{\delta \downarrow 0}
} \, \wh{\theta}_{\, \psi}(0,\delta) = 0$ as desired.
\end{proof}

In summary, starting from a nonifier $\psi(\bullet,\delta)$
satisfying the conditions in Lemma~\ref{lm:properties nonifier open}, the cumulative
distribution function $\wh{\theta}_{\, \psi}(\bullet,\delta)$ yields
a nondecreasing p-approximation function of the open Heaviside function satisfying conditions (A1), (A2), and (A3).

\gap

We next consider the reverse; i.e., we are given a bivariate function $\wh{\theta}(t,\delta)$ with $\wh{\theta}(\bullet,\delta)$ satisfying the conditions
in part (b) of
Proposition~\ref{pr:approx by truncation} and also (T2).  Since  $\wh{\theta}(\bullet,\delta)$ is B-differentiable, it is almost
everywhere differentiable.  Let $\psi(\bullet,\delta)$ be any integrable function such that $\psi(t,\delta) = \wh{\theta}(\bullet,\delta)^{\, \prime}(t)$
for almost all $t \in [ \, -\underline{\theta}(\delta),\overline{\theta}(\delta) \, ]$.  Defining $\psi(\bullet,\delta)$ to be zero outside the latter interval,
we deduce that
$\psi(\bullet,\delta)$ is a well-defined nonifier; moreover
\[
\wh{\theta}(t,\delta) \, = \, \displaystyle{
\int_{-\underline{\theta}(\delta)}^{\, t}
} \, \wh{\theta}(\bullet,\delta)^{\, \prime}(s) \, ds \, = \,
\displaystyle{
\int_{-\underline{\theta}(\delta)}^{\, t}
} \, \psi(s,\delta) \, ds \, = \, \wh{\theta}_{\psi}(t,\delta).
\]
Combining the two parts of the analysis, we conclude that the p-approximating functions for the open Heaviside function obtained
from the truncation approach coincide with those from the nonifier approach satisfying some mild properties.

\begin{example} \label{ex:truncated hinge} \rm
The truncated hinge loss function $T_{\rm h}(t,\delta/2) = \displaystyle{
\frac{1}{ \delta}
} \, \left[ \, \max( t + \delta/2, 0 ) - \max( t - \delta/2, 0 ) \, \right]$ is the averaged function
derived form the ``symmetric'' one-dimensional Steklov mollifier \cite[Definition 3.8]{ENWets95}:
$\psi(t,\delta) = \left\{ \begin{array}{ll}
1/\delta & \mbox{if $| \, t \, | \, \leq \, \delta/2$} \\ [3pt]
0 & \mbox{otherwise}
\end{array} \right.$, whose support is the interval $\left[ \, -\delta/2, \, \delta/2 \, \right]$.  As mentioned
before, the truncated hinge loss $T_{\rm h}(\bullet,\delta/2)$ fails to p-approximate the open Heaviside function because it violates condition (A1).
More interesting is the fact there do not exist nonifiers for this function that satisfy the conditions in Lemma~\ref{lm:properties nonifier open}.
Nevertheless, by considering the asymmetric Steklov function:
$\psi_a(t,\delta) = \left\{ \begin{array}{cl}
\displaystyle{
\frac{1}{\underline{\psi}(\delta) + \overline{\psi}(\delta)}
} & \mbox{if $t \, \in \, \left[ \, -\underline{\psi}(\delta), \, \overline{\psi}(\delta) \, \right]$} \\ [5pt]
0 & \mbox{otherwise}
\end{array} \right.$ with the end bounds $\underline{\psi}(\delta)$ and $\overline{\psi}(\delta)$ satisfying the
conditions in Lemma~\ref{lm:properties nonifier open},
we can derive a host of modified truncated hinge loss functions, such as (\ref{eq:nonifier for modified hinge}),
that p-approximate the open Heaviside function.  \hfill $\Box$
\end{example}

\section{Convergence to Pseudo B-stationary Solutions} \label{sec:approximation computation}

We consider the approximation of a pseudo B-stationary solution of the problem (\ref{eq:original POP})
by a combination of penalization of the functional constraint and approximation of the composite
Heaviside functions:
\begin{equation} \label{eq:approximated POP}
\displaystyle{
\operatornamewithlimits{\mbox{\bf minimize}}_{x \in X}
} \ \wh{\Phi}_{\lambda}(x,\delta) \, \triangleq \, c(x) + \underbrace{\displaystyle{
\sum_{k=1}^K
} \, \varphi_k(x) \, \theta_k^{\, \varphi}(g_k(x),\delta)}_{\mbox{denoted $\varphi(x,\delta)$}}
+ \lambda \, \max\left( \, \underbrace{\displaystyle{
\sum_{\ell=1}^L
} \, \phi_{\ell}(x) \, \theta_{\ell}^{\, \phi}(h_{\ell}(x),\delta) - b}_{\mbox{denoted $\phi(x,\delta)$}} , \, 0 \, \right),
\end{equation}
where for each pair $(k,\ell) \in [ K ] \times [ L ]$, $\theta_k^{\, \varphi}(\bullet,\delta)$ and
$\theta_{\ell}^{\, \phi}(\bullet,\delta)$ are p-approximations of the Heaviside function $\onebld_{( \, 0,\infty \, )}$
with support in the interval $[ \, -\underline{\theta}_k^{\, \varphi}(\delta),\overline{\theta}_k^{\, \varphi}(\delta) \, ]$
and $[ \, -\underline{\theta}_{\ell}^{\, \phi}(\delta),\overline{\theta}_{\ell}^{\, \phi}(\delta) \, ]$, respectively,
that shrink to zero when $\delta \downarrow 0$;
in particular, (A1), (A2), and (A3) in Definition~\ref{def:approximation} are satisfied by these approximating functions.
With the functions $c$ and each $\varphi_k$ being Lipschitz continuous on $X$ with Lipschtiz constants $\mbox{Lip}_c$ and
$\mbox{Lip}_{\varphi}$, respectively, we take $\lambda$ satisfying (\ref{eq:condition on penalty}).

\gap

One of the
conditions that we will impose in the analysis is a sign condition, labelled (C2) below, on the directional derivatives of the
functions $g_k$ and $h_{\ell}$.  The lemma below shows that this will hold if these are convex piecewise affine functions.

\begin{lemma} \label{lm:zero dd of convex/concave} \rm
Let $f : \mathbb{R}^n \to \mathbb{R}$ be a convex piecewise affine function.  For every vector $\bar{x}$, there exists a neighborhood ${\cal N}$
of $\bar{x}$ such that
for all $v \in \mathbb{R}^n$,
\[
f^{\, \prime}(\bar{x};v) \, \leq \, 0 \ \Rightarrow \ f^{\, \prime}(x;v) \, \leq \, 0, \epc \forall \, x \, \in \, {\cal N}.
\]
\end{lemma}

\begin{proof}  Since $f$ is convex piecewise affine, it can be written as the pointwise maximum of finitely many affine functions
\cite[Proposition~4.4.6]{CuiPang2021}; i.e.,
\[
f(x) \, = \, \displaystyle{
\max_{1 \leq i \leq I}
} \, \left( ( a^i )^{\top} x + b_i \right), \epc \forall \, x \, \in \, \mathbb{R}^n,
\]
for some positive integer $I$, $n$-vectors $\{ a^{\, i} \}_{i=1}^I$, and scalars $\{ b_i \}_{i=1}^I$.  For any $x$, let
\[
{\cal A}(x) \, \triangleq \, \left\{ \, i \, \in \, [ I ] \, \mid \, f(x) \, = \, ( a^i )^{\top} x + b_i \, \right\}
\]
be the maximizing index set of the affine pieces of $f$.  It then follows that for the given $\bar{x}$, there exists a neighborhood ${\cal N}$ such that
${\cal A}(x) \subseteq {\cal A}(\bar{x})$ for all $x \in {\cal N}$.  Since
\[
f^{\, \prime}(x;v) \, = \, \displaystyle{
\max_{i \in {\cal A}(x)}
} \, \left( ( a^i )^{\top} v + b_i \right), \epc \forall \, ( x,v ) \, \in \, \mathbb{R}^n \times \mathbb{R}^n,
\]
the desired conclusion of the lemma follows readily.
\end{proof}

Let $\{ \delta_{\nu} \}$ be a sequence of positive scalars converging to zero.
Let $\{ x^{ \nu} \}$ be a corresponding sequence of d-stationary points
of $\wh{\Phi}_{\lambda}(\bullet,\delta_{\nu})$ on $X$ so that $\wh{\Phi}_{\lambda}(\bullet,\delta_{\nu})^{\, \prime}(x^{\nu};v) \geq 0$ for all
$v \in {\cal T}(X;x^{\nu})$.
Suppose that $\{ x^{\, \nu} \}$ converges to the limit $x^*$.  In what follows, we show that $x^*$ is a pseudo B-stationary point
of (\ref{eq:original POP}) by verifying two things under the condition (\ref{eq:condition on penalty}) on $\lambda$:
(a) $x^*$ is feasible to (\ref{eq:original POP}), and (b) the implication (\ref{eq:ACQ stationarity}), which we
restate below in terms of the vector $x^*$ on hand:
%
\begin{equation} \label{eq:barx stationary approx}
v \, \in {\cal L}(S_{\rm ps}(x^*);x^*) \ \Rightarrow \ \Phi(\bullet;x^*)^{\, \prime}(x^*;v) \, \geq \, 0.
\end{equation}
Letting $\sigma_k(x;v) \triangleq \mbox{sgn}( g_k^{\, \prime}(x;v) )$
with $\mbox{sgn}(0)$ defined to be zero,
we write the directional derivative of $\varphi(\bullet,\delta_{\nu})$ at $x^{\nu}$ as the sum of 3 terms:
\[ \begin{array}{l}
\varphi(\bullet,\delta_{\nu})^{\, \prime}(x^{\, \nu};v) \\ [0.1in]
= \, \displaystyle{
\sum_{k=1}^K
} \, \left[ \, \varphi_k^{\, \prime}(x^{\, \nu};v) \, \theta_k^{\, \varphi}(g_k(x^{\, \nu}),\delta_{\nu}) +
\varphi_k(x^{\, \nu}) \,
| \, g_k^{\, \prime}(x^{\, \nu};v) \, | \,
( \theta_k^{\, \varphi}(\bullet,\delta_{\nu}) )^{\prime}( g_k(x^{\, \nu});\sigma_k(x^{\, \nu};v) ) \, \right] \\ [0.25in]
= \, T_>^{\, \varphi}(x^{\, \nu};v) + T_=^{\, \varphi}(x^{\, \nu};v) + T_<^{\, \varphi}(x^{\, \nu};v), \epc \mbox{where}
\end{array} \]
\[ \begin{array}{l}
T_>^{\, \varphi}(x^{\, \nu};v) \, = \,
\displaystyle{
\sum_{k \in {\cal K}_>(x^*)}
} \, \left[ \, \varphi_k^{\, \prime}(x^{\, \nu};v) \, \theta_k^{\, \varphi}(g_k(x^{\, \nu}),\delta_{\nu}) + \varphi_k(x^{\, \nu}) \,
| \, g_k^{\, \prime}(x^{\, \nu};v) \, | \,
( \theta_k^{\, \varphi}(\bullet,\delta_{\nu}) )^{\, \prime}( g_k(x^{\, \nu});\sigma_k(x^{\, \nu};v) ) \, \right] \\ [0.25in]
T_=^{\, \varphi}(x^{\, \nu};v) \, = \,
\displaystyle{
\sum_{k \in {\cal K}_=(x^*)}
} \, \left[ \, \varphi_k^{\, \prime}(x^{\, \nu};v) \, \theta_k^{\, \varphi}(g_k(x^{\, \nu}),\delta_{\nu}) + \varphi_k(x^{\, \nu}) \,
| \, g_k^{\, \prime}(x^{\, \nu};v) \, | \,
( \theta_k^{\, \varphi}(\bullet,\delta_{\nu}) )^{\, \prime}( g_k(x^{\, \nu});\sigma_k(x^{\, \nu};v) ) \, \right] \\ [0.25in]
T_<^{\, \varphi}(x^{\, \nu};v) \, = \,
\displaystyle{
\sum_{k \in {\cal K}_<(x^*)}
} \, \left[ \, \varphi_k^{\, \prime}(x^{\, \nu};v) \, \theta_k^{\, \varphi}(g_k(x^{\, \nu}),\delta_{\nu}) + \varphi_k(x^{\, \nu}) \,
| \, g_k^{\, \prime}(x^{\, \nu};v) \, | \,
( \theta_k^{\varphi}(\bullet,\delta_{\nu}))^{\, \prime}( g_k(x^{\, \nu});\sigma_k(x^{\, \nu};v) ) \, \right].
\end{array}
\]
Consider the sum $T_>^{\varphi}(x^{\, \nu};v)$.  For each $k \in {\cal K}_>(x^*)$, we have $g_k(x^*) > 0$.
Hence $g_k(x^{\, \nu}) > \overline{\theta}_k^{\, \varphi}(\delta_{\nu})$ for all $\nu$ sufficiently large;
this yields, by (A1):
\[
\theta_k(g_k(x^{\, \nu}),\delta_{\nu}) \, = \, 1 \ \mbox{ and }
\theta_k^{\, \prime}(\bullet,\delta_{\nu})^{\, \prime}( g_k(x^{\, \nu});\sigma_k(x^{\, \nu};v) ) \, = \, 0.
\]
Thus,
\[
T_>^{\varphi}(x^{\, \nu};v) \, = \, \displaystyle{
\sum_{k \in {\cal K}_>(x^*)}
} \, \phi_k^{\, \prime}(x^{\, \nu};v).
\]
Similarly, we can show that $T_<^{\varphi}(x^{\, \nu};v) = 0$ for all $\nu$ sufficiently large.  The rest of the proof is
divided into three parts.

\gap

$\bullet $ To analyze the term $T_=^{\varphi}(x^{\, \nu};v)$ and the corresponding term in $\phi(\bullet,\delta_{\nu})^{\, \prime}(x^*;v)$,
we need assumptions (C1) and (C2):

\gap

(C1): There exists an open neighborhood ${\cal N}_*$ such that on $X \cap {\cal N}_*$, the functions
$\{ \varphi_k \}_{k \in {\cal K}_=(x^*)}$ and $\left\{ \, \phi_{\ell} \, \right\}_{\ell \in {\cal L}_=(x^*)}$
are  nonnegative; this is a pointwise sign condition related to those in
Proposition~\ref{pr:lsc case}, demanding in particular the nonnegativity of these functions in
a region around $x^*$ and possibly outside the respective sets
$X \cap g_k^{-1}( \, -\infty, 0 \, ]$ and $X \cap h_{\ell}^{-1}( \, \infty, 0 \, ]$.

\gap

(C2): The implications
\begin{equation} \label{eq:strong convex gk}
\left. \begin{array}{l}
g_k^{\, \prime}(x^*;v) \, \leq \, 0 \epc \forall \, k \in {\cal K}_=(x^*) \\ [5pt]
v \, \in \, {\cal T}(X;x^*)
\end{array} \right\} \ \Rightarrow \ g_k^{\, \prime}(x;v) \, \leq \, 0 \epc \mbox{for all $k \in {\cal K}_=(x^*)$ and all $x \in {\cal N}_*$}
\end{equation}
\begin{equation} \label{eq:strong convex hell}
\left. \begin{array}{l}
h_{\ell}^{\, \prime}(x^*;v) \, \leq \, 0 \epc \forall \, \ell \in {\cal L}_=(x^*) \\ [5pt]
v \, \in \, {\cal T}(X;x^*)
\end{array} \right\} \ \Rightarrow \ h_{\ell}^{\, \prime}(x;v) \, \leq \, 0 \epc \mbox{for all $\ell \in {\cal L}_=(x^*)$ and all $x \in {\cal N}_*$}
\end{equation}
are motivated by
Lemma~\ref{lm:zero dd of convex/concave} which provides sufficient conditions for them to hold.

\gap

Under (C1) and (C2), we deduce that if $v$ satisfies the left-hand condition in (\ref{eq:strong convex gk}),
then by (\ref{eq:A3 for p-approx}) applied to $\theta_k^{\, \varphi}$ , we obtain, for all $\nu$ sufficiently large
and all $k \in {\cal K}_=(x^*)$,
\[
\varphi_k(x^{\, \nu}) \, | \, g_k^{\, \prime}(x^{\, \nu};v) \, | \,
\theta_k^{\, \varphi}(\bullet,\delta_{\nu})^{\, \prime}( g_k(x^{\, \nu});\sigma_k(x^{\, \nu};v) ) \, \leq \, 0
\]
Hence, for such $\nu$, we have
\[
T_{k;=}^{\varphi}(x^{\, \nu};v) \, \leq \,
\varphi_k^{\, \prime}(x^{\, \nu};v) \, \theta_k^{\, \varphi}(g_k(x^{\, \nu}),\delta_{\nu}).
\]
Consequently, we deduce, for all $\nu $ sufficiently large,
\begin{equation} \label{eq:varphi dd ub}
\begin{array}{l}
\left. \begin{array}{l}
g_k^{\, \prime}(x^*;v) \, \leq \, 0 \epc \forall \, k \in {\cal K}_=(x^*) \\ [5pt]
v \, \in \, {\cal T}(X;x^*)
\end{array} \right\} \\ [0.2in]
\hspace{0.3in} \Rightarrow \ \varphi(\bullet,\delta_{\nu})^{\, \prime}(x^{\, \nu},v)
\, \leq \, \displaystyle{
\sum_{k \in {\cal K}_>(x^*)}
} \, \varphi_k^{\, \prime}(x^{\, \nu};v) + \displaystyle{
\sum_{k \in{\cal K}_=(x^*)}
} \, \varphi_k^{\, \prime}(x^{\, \nu};v) \, \theta_k^{\, \varphi}(g_k(x^{\, \nu}),\delta_{\nu}).
\end{array} \end{equation}
Similarly, we also have, for all $\nu$ sufficiently large,
\begin{equation} \label{eq:phi dd ub}
\begin{array}{l}
\left. \begin{array}{l}
h_{\ell}^{\, \prime}(x^*;v) \, \leq \, 0 \epc \forall \, \ell \in {\cal L}_=(x^*) \\ [5pt]
v \, \in \, {\cal T}(X;x^*)
\end{array} \right\} \\ [0.2in]
\hspace{0.3in} \Rightarrow \ \phi(\bullet,\delta_{\nu})^{\, \prime}(x^{\, \nu},v)
\, \leq \, \displaystyle{
\sum_{\ell \in {\cal L}_>(x^*)}
} \, \phi_{\ell}^{\, \prime}(x^{\, \nu};v) + \displaystyle{
\sum_{\ell \in {\cal L}_=(x^*)}
} \, \phi_{\ell}^{\, \prime}(x^{\, \nu};v) \, \theta_k^{\, \phi}(h_{\ell}(x^{\, \nu}),\delta_{\nu}).
\end{array} \end{equation}
$\bullet $ To establish the feasibility of $x^*$ for the problem (\ref{eq:original POP}), we postulate three more assumptions:
as we will see below, assumption (C3) ensures the objective recovery of the sequence $\{ x^{\, \nu} \}$; (C4) is
the same as that in Theorem~\ref{th:fixed-point stationary} and (C5) is the Clarke regularity on several key functions.

\gap

(C3) For all $k \in {\cal K}_=(x^*)$ and $\ell \in {\cal L}_=(x^*)$, the limits
\begin{equation} \label{eq:no gap}
\displaystyle{
\lim_{\nu \to \infty}
} \, \theta_k^{\, \varphi}(g_k(x^{\, \nu}),\delta_{\nu}) \, = \, 0 \epc \mbox{and} \epc
\displaystyle{
\lim_{\nu \to \infty}
} \, \theta_{\ell}^{\, \phi}(h_{\ell}(x^{\, \nu}),\delta_{\nu}) \, = \, 0
\end{equation}
are {\sl functional consistency} requirements of the sequence $\{ x^{\, \nu} \}$ in the following sense.
These limits ensure that for all pairs $( k,\ell)$,
\[
\displaystyle{
\lim_{\nu \to \infty}
} \, \theta_k^{\, \varphi}(g_k(x^{\, \nu}),\delta_{\nu}) = \onebld_{( \, 0,\infty \, )}(g_k(x^*)) \epc \mbox{and} \epc
\displaystyle{
\lim_{\nu \to \infty}
} \, \theta_{\ell}^{\, \phi}(h_{\ell}(x^{\, \nu}),\delta_{\nu}) = \onebld_{( \, 0,\infty \, )}(h_{\ell}(x^*)).
\]
Together with the same limits for $k \in {\cal K}_>(x^*) \cup {\cal K}_<(x^*)$ and $\ell \in {\cal L}_>(x^*) \cup {\cal L}_<(x^*)$,
we deduce in particular that
\[
\displaystyle{
\lim_{\nu \to \infty}
} \, \left[ \, c(x^{\, \nu}) + \displaystyle{
\sum_{k=1}^K
} \, \varphi_k(x^*) \, \theta_k^{\, \varphi}(g_k(x^{\, \nu}),\delta_{\nu}) \, \right] \, = \, \Phi(x^*),
\]
which is reasonable to postulate in order for
$x^*$ to be a stationarity point of some kind for the objective function of the original problem (\ref{eq:original POP}).
If the approximation function $\theta_k$ is chosen as in part (c) of Proposition~\ref{pr:approx by truncation}, i.e. if
\[
\theta_k(t,\delta) \, \triangleq \, T_{[ \, 0,1 \, ]}\left( \psi_k\left( q_k(\delta) + \displaystyle{
\frac{t}{m_k(\delta)}
} \, \right) \, \right)
\]
where the functions $\psi_k$, $q_k$, and $m_k$ are as specified in the proposition, then the limit (\ref{eq:no gap})
holds if $\displaystyle{
\lim_{\nu \to \infty}
} \, \displaystyle{
\frac{g_k(x^{\nu})}{m_k(\delta_{\nu})}
} \, = \, 0$.  Admittedly, this is a condition on the sequence of iterates $\{ x^{\nu} \}$ relative to the sequence of
parameters $\{ \delta_{\nu} \}$ that needs to be addressed from the source of the iterates, e.g., as generated by an iterative algorithm.
Regrettably, the design of such algorithms is beyond the scope of the present work but will be the focus of a subsequent
computational study which will be guided by this high-level background result.  The limit on the
sequence $\{ \theta_{\ell}^{\, \phi}(h_{\ell}(x^{\, \nu}),\delta_{\nu}) \}$ has to do with the satisfaction of the functional constraint
and facilitates the demonstration of the desired pseudo B-stationarity of the limit $x^*$; this
will become clear in what follows.  Subsequently, we will establish a weaker stationarity property of the limit $x^*$ without (C3).

\gap

 {(C4) For any $x\in X$ such that $\displaystyle{
\sum_{\ell \in {\cal L}_>(x)}
} \, \phi_{\ell}(x)> b$,
there exists a vector
$\bar{v} \in {\cal T}( \wh{S}_{\rm ps}(x);x)$ with unit length satisfying:
$\displaystyle{
\sum_{\ell \in {\cal L}_>(x)}
} \, \phi_{\ell}^{\, \prime}(x;\bar{v}) \leq -1$;} and

\gap

(C5) the functions $c$, $\{ \varphi_k \}_{k \in {\cal K}_>(x^*)}$, and $\{ \phi_{\ell} \}_{\ell \in {\cal L}_>(x^*)}$
are Clarke regular \cite{Clarke83} at $x^*$.  In the context of a B-differentiable
function $f : {\cal O} \to \mathbb{R}$ (which these functions are), Clarke regularity at a vector $\bar{x} \in {\cal O}$
means that for all sequences $\{ z^{\, \nu} \}$ converging to $\bar{x}$, it holds that
\[
\displaystyle{
\limsup_{\nu \to \infty}
} \,  f^{\, \prime}(z^{\, \nu};v) \, \leq \, f^{\, \prime}(\bar{x};v), \epc \forall \, v \, \in \, \mathbb{R}^n.
\]
 {
Condition (C4) requires that for any $x\in X$ that is infeasible to problem (\ref{eq:original POP}), the constraint function $\displaystyle{
\sum_{\ell\in {\cal L}_>(x)}
} \, \phi_{\ell}(\bullet)$
has a descent direction at $x$.  This condition is in the same spirit as a classical one since the early days
of the theory of exact penalty methods \cite{PilloFacchinei89}
for the recovery of feasibility (and hence stationarity) in a penalized problem; see also \cite[Proposition~9.2.2 (a)]{CuiPang2021}.
As in problems without the Heaviside function, condition (C4) can be related to the concept of weak-sharp minima \cite{Ferris88}
properly extended.  We give a brief discussion of the connection in an Appendix.
The Clarke regularity in (C5) holds in particular if the considered functions are the composition of a convex function
with a smooth mapping \cite[Exercise~10.25]{RockafellarWets09}.}


Continuing the analysis, assume for the sake of contradiction that $\displaystyle{
\sum_{\ell=1}^L
} \, \phi_{\ell}(x^*) \, \onebld_{( \, 0,\infty )}(h_{\ell}(x^*)) > b$.  Then for all $\nu$ sufficiently large,
$\phi(x^{\, \nu},\delta_{\nu}) \, = \displaystyle{
\sum_{\ell=1}^L
} \, \phi_{\ell}(x^{\, \nu}) \, \theta_{\ell}^{\, \varphi}(h_{\ell}(x^{\, \nu}),\delta_{\nu}) - b > 0$.
Hence,
for the vector $\bar{v}$ in (C4), we have
\begin{equation} \label{eq:consequence of C4}
\begin{array}{lll}
0 & \leq & \wh{\Phi}_{\lambda}(\bullet,\delta_{\nu})^{\, \prime}(x^{\, \nu};\bar{v}) \, = \, c^{\, \prime}(x^{\, \nu};\bar{v}) +
\varphi(\bullet,\delta_{\nu})^{\, \prime}(x^{\, \nu};\bar{v}) + \lambda \, \phi(\bullet,\delta_{\nu})^{\, \prime}(x^{\, \nu};\bar{v})
\\ [0.1in]
& \leq & c^{\, \prime}(x^{\, \nu};\bar{v}) +  \displaystyle{
\sum_{k \in {\cal K}_>(x^*)}
} \, \varphi_k^{\, \prime}(x^{\, \nu};\bar{v}) + \displaystyle{
\sum_{k \in{\cal K}_=(x^*)}
} \, \varphi_k^{\, \prime}(x^{\, \nu};\bar{v}) \, \theta_k^{\, \varphi}(g_k(x^{\, \nu}),\delta_{\nu}) \\ [0.3in]
& & \hspace{0.7in} + \, \lambda \, \left[ \,
\displaystyle{
\sum_{\ell \in {\cal L}_>(x^*)}
} \, \phi_{\ell}^{\, \prime}(x^{\, \nu};\bar{v}) + \displaystyle{
\sum_{\ell \in {\cal L}_=(x^*)}
} \, \phi_{\ell}^{\, \prime}(x^{\, \nu};\bar{v}) \, \theta_k^{\, \phi} (h_{\ell}(x^{\, \nu}),\delta_{\nu}) \, \right] \\ [0.3in]
& \leq & ( \, \mbox{Lip}_c + K \, \mbox{Lip}_{\varphi} \, ) +
\lambda \, \left[ \,
\displaystyle{
\sum_{\ell \in {\cal L}_>(x^*)}
} \, \phi_{\ell}^{\, \prime}(x^{\, \nu};\bar{v}) + \mbox{Lip}_{\phi} \, \displaystyle{
\sum_{\ell \in {\cal L}_=(x^*)}
} \, \theta_k^{\, \phi}(h_{\ell}(x^{\, \nu}),\delta_{\nu}) \, \right],
\end{array} \end{equation}
where $\mbox{Lip}_{\phi}$ is a local Lipschitz constant of $\phi_{\ell}$ near $x^*$.  Taking the limit $\nu \to \infty$ and using
(\ref{eq:no gap}), we obtain, by the Clarke regularity of the functions $\{ \phi_{\ell} \}_{\ell \in {\cal L}_>(x^*)}$ at $x^*$,
\[ \begin{array}{lll}
0 & \leq & \mbox{Lip}_c + K \, \mbox{Lip}_{\varphi} + \lambda \, \displaystyle{
\sum_{\ell \in {\cal L}_>(x^*)}
} \, \displaystyle{
\limsup_{\nu \to \infty}
} \, \phi_{\ell}^{\, \prime}(x^{\, \nu};\bar{v}) \\ [0.2in]
& \leq &  \mbox{Lip}_c + K \, \mbox{Lip}_{\varphi} + \lambda \, \displaystyle{
\sum_{\ell \in {\cal L}_>(x^*)}
} \, \phi_{\ell}^{\, \prime}(x^*;\bar{v}) \, \leq \,  \mbox{Lip}_c + K \, \mbox{Lip}_{\varphi} - \lambda,
\end{array} \]
This contradiction completes the feasibility proof of $x^*$ for (\ref{eq:original POP}).

\gap

$\bullet $ Finally, to complete the proof of the desired implication (\ref{eq:barx stationary approx}), suppose
that $x^*$ satisfies the functional constraint as an equality; i.e.,
$\displaystyle{
\sum_{\ell=1}^L
} \, \phi_{\ell}(x^*) \, \onebld_{( \, 0,\infty )}(h_{\ell}(x^*)) = b$.
Let $v \in {\cal L}(S_{\rm ps}(x^*);x^*)$ be arbitrary.  We have
\[ \begin{array}{l}
\displaystyle{
\limsup_{\nu \to \infty}
} \, \left[ \, \max\left( \phi(\bullet,\delta_{\nu}), \, 0 \, \right) \, \right]^{\, \prime}(x^{\nu};v) \, \leq \, \displaystyle{
\limsup_{\nu \to \infty}
} \, \max\left( \phi(\bullet,\delta_{\nu})^{\, \prime}(x^{\, \nu};v), \, 0 \, \right), \epc \mbox{by Lemma~\ref{lm:dd of max composite}} \\ [0.1in]
\epc \leq \, \displaystyle{
\limsup_{\nu \to \infty}
} \, \max\left( \, \displaystyle{
\sum_{\ell \in {\cal L}_>(x^*)}
} \, \phi_{\ell}^{\, \prime}(x^{\, \nu};v) + \displaystyle{
\sum_{\ell \in {\cal L}_=(x^*)}
} \, \phi_{\ell}^{\, \prime}(x^{\, \nu};v) \, \theta_{\ell}^{\, \phi}(h_{\ell}(x^{\, \nu}),\delta_{\nu}), \, 0 \, \right) \epc \mbox{by
(\ref{eq:phi dd ub})} \\ [0.3in]
\epc \leq \, \max\left( \, \displaystyle{
\limsup_{\nu \to \infty}
} \, \left[ \, \displaystyle{
\sum_{\ell \in {\cal L}_>(x^*)}
} \, \phi_{\ell}^{\, \prime}(x^{\, \nu};v) + \mbox{Lip}_{\phi} \, \displaystyle{
\sum_{\ell \in {\cal L}_=(x^*)}
} \, \theta_{\ell}^{\, \phi}(h_{\ell}(x^{\, \nu}),\delta_{\nu}) \, \right], \ 0 \, \right) \\ [0.35in]
\epc \leq \, \max\left( \,  \displaystyle{
\sum_{\ell \in {\cal L}_>(x^*)}
} \, \phi_{\ell}^{\, \prime}(x^*;v), \, 0 \, \right), \epc \mbox{by Clarke regularity and (C3)} \\ [0.3in]
\epc = \, 0, \epc \mbox{because $\displaystyle{
\sum_{\ell \in {\cal L}_>(x^*)}
} \, \phi_{\ell}^{\, \prime}(x^*;v) \leq 0$ is a stipulation on $v$.}
\end{array} \]
In the other case where $\displaystyle{
\sum_{\ell=1}^L
} \, \phi_{\ell}(x^*) \, \onebld_{( \, 0,\infty )}(h_{\ell}(x^*)) < b$, we must have $\phi(x^{\, \nu},\delta_{\nu}) < b$
for all $\nu$ sufficiently large.  Thus, $\left[ \, \max\left( \phi(\bullet,\delta_{\nu}), \, 0 \, \right) \, \right]^{\, \prime}(x^{\nu};v) = 0$
for all such $\nu$.  Consequently, in either case, it follows that
\[ \begin{array}{lll}
0 & \leq & \displaystyle{
\limsup_{\nu \to \infty}
} \, \left[ \, c^{\, \prime}(x^{\, \nu};v) + \varphi(\bullet,\delta_{\nu})^{\, \prime}(x^{\, \nu};v) \, \right] \\ [0.1in]
& \leq & c^{\, \prime}(x^*;v) + \displaystyle{
\sum_{k \in {\cal K}_>(x^*)}
} \, \varphi_k^{\, \prime}(x^*;v), \epc \mbox{by (\ref{eq:varphi dd ub}), (C3), and Clarke regularity},
\end{array} \]
establishing the right-hand side of (\ref{eq:barx stationary approx}).  We have thus proved the following
main result of this section.

\begin{theorem} \label{th:ccp-stationary} \rm
Under the blanket assumption of problem (\ref{eq:original POP}),
let $c$ and each $\varphi_k$ be Lipschitz continuous on $X$ with Lipschitz constants $\mbox{Lip}_c$ and
$\mbox{Lip}_{\varphi}$, respectively.
Let $\{ \delta_{\nu} \}$ be a sequence of positive scalars converging to zero and for each $\nu$, let
$x^{\nu}$ be a d-stationary point of $\wh{\Phi}_{\lambda}(\bullet,\delta_{\nu})$ on $X$ with $\lambda$ satisfying
(\ref{eq:condition on penalty}).  Suppose $\displaystyle{
\lim_{\nu \to \infty}
} \, x^{\, \nu} = x^*$.  Under assumptions (C1)--(C5), it holds that $x^*$ is a pseudo B-stationary solution of (\ref{eq:original POP}).
\hfill $\Box$
\end{theorem}

Without (C3), the desired pseudo B-stationarity of the limit $x^*$ as defined in Definition~\ref{df:pseudo concepts} is in jeopardy.
By strengthening (C4) and (C5), we can still establish a weak pseudo B-stationarity property of $x^*$.  The two strengthened conditions are:

\gap

(C4$^{\, \prime}$) There exists a vector
$\bar{v} \in {\cal T}( \wh{S}_{\rm ps}(x^*);x^*)$ with unit length satisfying:
\begin{equation} \label{eq:strengthened C4}
\displaystyle{
\sum_{\ell \in {\cal L}_>(\bar{x})}
} \, \phi_{\ell}^{\, \prime}(x^*;\bar{v}) + \displaystyle{
\sum_{\ell \in {\cal L}_=(\bar{x})}
} \, \max\left( \, \phi_{\ell}^{\, \prime}(x^*;\bar{v}),0 \, \right) \, \leq \, -1.
\end{equation}
While more demanding than (C4), this strengthened condition is still in the spirit of the common requirements in the theory of exact penalization
(see conditions (a) and (b) in \cite[Theorem~9.2.1]{CuiPang2021}).

\gap

(C5$^{\, \prime}$) In addition to those in (C5), the functions $\{ \varphi_k \}_{k \in {\cal K}_=(x^*)}$, and $\{ \phi_{\ell} \}_{\ell \in {\cal L}_=(x^*)}$
are also Clarke regular at $x^*$. 

\begin{proposition} \label{pr:weak pseudo Bstationary} \rm
Under the blanket assumption of problem (\ref{eq:original POP}),
let $c$ and each $\varphi_k$ be Lipschitz continuous on $X$ with Lipschitz constants $\mbox{Lip}_c$ and
$\mbox{Lip}_{\varphi}$, respectively.
Let $\{ \delta_{\nu} \}$ be a sequence of positive scalars converging to zero and for each $\nu$, let
$x^{\nu}$ be a d-stationary point of $\wh{\Phi}_{\lambda}(\bullet,\delta_{\nu})$ on $X$ with $\lambda$ satisfying
(\ref{eq:condition on penalty}).  Suppose $\displaystyle{
\lim_{\nu \to \infty}
} \, x^{\, \nu} = x^*$.  Under assumptions (C1), (C2), (C4$^{\, \prime}$), and (C5$^{\, \prime}$), $x^*$ is feasible to
(\ref{eq:original POP}) and there exist scalars
$\{ \, \xi_k^* \, \}_{k \in {\cal K}_=(x^*)} \, \cup \, \{ \, \mu_{\ell}^* \, \}_{\ell \in {\cal L}_=(x^*)} \, \subset \, [ \, 0, 1 \, ]$
such that $x^*$ is a B-stationary solution of:
\begin{equation} \label{eq:weak barx stationary problem}
\begin{array}{l}
\displaystyle{
\operatornamewithlimits{\mbox{\bf minimize}}_{x \in X}
} \  \Phi_{\geq}^{\xi^*}(x;x^*) \, \triangleq \, c(x) + \displaystyle{
\sum_{k \, \in \, {\cal K}_>(x^*)}
} \, \varphi_k(x) + \underbrace{\displaystyle{
\sum_{k \, \in \, {\cal K}_=(x^*)}
} \, \xi_k^* \, \varphi_k(x)}_{\mbox{extra term}} \\ [0.5in]
\left. \begin{array}{ll}
\mbox{\bf subject to} & \displaystyle{
\sum_{\ell \, \in \, {\cal L}_>(x^*)}
} \, \phi_{\ell}(x) + \underbrace{\displaystyle{
\sum_{\ell \, \in \, {\cal L}_=(x^*)}
} \, \mu_{\ell}^* \, \phi_{\ell}(x)}_{\mbox{extra term}} \, \leq \, b \\ [0.5in]
& g_k(x) \, \leq \, 0 \epc \forall \, k \in {\cal K}_{\leq}(x^*)
\\ [0.1in]
& g_k(x) \, \geq \, 0 \epc \forall \, k \in {\cal K}_>(x^*) \\ [0.1in]
& h_{\ell}(x) \, \leq \, 0 \epc \forall \, \ell \in {\cal L}_{\leq}(x^*) \\ [0.1in]
\mbox{\bf and} & h_{\ell}(x) \, \geq \, 0 \epc \forall \, \ell \in {\cal L}_>(x^*)
\end{array} \right\} .
\end{array} \end{equation}
\end{proposition}

\begin{proof}  Belonging to the interval $[ \, 0,1 \, ]$, the sequences
\[
\left\{ \, \theta_k^{\, \varphi} (g_k(x^{\, \nu}),\delta_{\nu}) \, \right\}_{k \in {\cal K}_=(x^*)} \epc \mbox{and} \epc
\left\{ \, \theta_{\ell}^{\, \phi} (h_{\ell}(x^{\, \nu}),\delta_{\nu}) \, \right\}_{\ell \in {\cal L}_=(x^*)}
\]
have accumulation points, say $\{ \xi_k^* \}_{k \in {\cal K}_=(x^*)}$ and $\{ \mu_{\ell}^* \}_{\ell \in {\cal L}_=(x^*)}$,
which we may assume, without loss of generality are limits of the displayed sequences, respectively.
We claim that $x^*$ satisfies the functional constraint in (\ref{eq:weak barx stationary problem}).
Assume otherwise.  Then $\displaystyle{
\sum_{\ell \in {\cal L}_>(x^*)}
} \, \phi_{\ell}(x^*) + \displaystyle{
\sum_{\ell \in {\cal L}_=(x^*)}
} \, \mu^* \, \phi_{\ell}(x^*) > b$.  Since
\[
\displaystyle{
\lim_{\nu \to \infty}
} \, \phi(x^{\, \nu},\delta_{\nu}) \, = \, \displaystyle{
\sum_{\ell \in {\cal L}_>(x^*)}
} \, \phi_{\ell}(x^*) + \displaystyle{
\sum_{\ell \in {\cal L}_=(x^*)}
} \, \mu^* \, \phi_{\ell}(x^*)
\]
it follows that $\phi(x^{\, \nu},\delta_{\nu}) > b$ for all $\nu$ sufficiently large.
Continuing from (\ref{eq:consequence of C4}), we have
\[
\begin{array}{lll}
0 & \leq & \wh{\Phi}_{\lambda}(\bullet,\delta_{\nu})^{\, \prime}(x^{\, \nu};\bar{v})
\, \leq \, c^{\, \prime}(x^{\, \nu};\bar{v}) +  \displaystyle{
\sum_{k \in {\cal K}_>(x^*)}
} \, \varphi_k^{\, \prime}(x^{\, \nu};\bar{v}) + \displaystyle{
\sum_{k \in{\cal K}_=(x^*)}
} \, \varphi_k^{\, \prime}(x^{\, \nu};\bar{v}) \, \theta_k^{\, \varphi}(g_k(x^{\, \nu}),\delta_{\nu}) \\ [0.3in]
& & \hspace{1.5in} + \, \lambda \, \left[ \,
\displaystyle{
\sum_{\ell \in {\cal L}_>(x^*)}
} \, \phi_{\ell}^{\, \prime}(x^{\, \nu};\bar{v}) + \displaystyle{
\sum_{\ell \in {\cal L}_=(x^*)}
} \, \phi_{\ell}^{\, \prime}(x^{\, \nu};\bar{v}) \, \theta_{\ell}^{\, \phi} (h_{\ell}(x^{\, \nu}),\delta_{\nu}) \, \right] \\ [0.3in]
& \leq & ( \, \mbox{Lip}_c + K \, \mbox{Lip}_{\varphi} \, ) +
\lambda \, \left[ \,
\displaystyle{
\sum_{\ell \in {\cal L}_>(x^*)}
} \, \phi_{\ell}^{\, \prime}(x^{\, \nu};\bar{v}) + \displaystyle{
\sum_{\ell \in {\cal L}_=(x^*)}
} \, \max\left( \, \phi_{\ell}^{\, \prime}(x^{\, \nu};\bar{v}), \, 0 \, \right) \, \right] \\ [0.3in]
& & \epc \mbox{because both $\theta_k^{\, \varphi}(g_k(x^{\, \nu}),\delta_{\nu})$ and $\theta_{\ell}^{\, \phi} (h_{\ell}(x^{\, \nu}),\delta_{\nu})$
are in $[ \, 0,1 \, ]$}.
\end{array} \]
The strengthened condition (\ref{eq:strengthened C4}) and the Clarke regularity condition (C5$^{\, \prime}$) then yield a contradiction by
letting $\nu \to \infty$.  This shows that $\displaystyle{
\sum_{\ell \in {\cal L}_>(x^*)}
} \, \phi_{\ell}(x^*) + \displaystyle{
\sum_{\ell \in {\cal L}_=(x^*)}
} \, \mu^*_{\ell} \, \phi_{\ell}(x^*) \leq b$.  Since $\mu^*_{\ell} \, \phi_{\ell}(x^*) \geq 0$ for all $\ell \in {\cal L}_>(x^*)$, we have
\[
\displaystyle{
\sum_{\ell \in {\cal L}_>(x^*)}
} \, \phi_{\ell}(x^*) + \displaystyle{
\sum_{\ell \in {\cal L}_=(x^*)}
} \, \mu^*_{\ell} \, \phi_{\ell}(x^*) \, \geq \, \displaystyle{
\sum_{\ell=1}^L
} \, \phi_{\ell}(x^*) \, \onebld_{( \, 0,\infty \, )}(h_{\ell}(x^*)),
\]
thus $x^*$ is feasible to (\ref{eq:original POP}).  Finally, the proof that $x^*$ is a B-stationary solution of
(\ref{eq:weak barx stationary problem}) is similar to that in the last part of the proof of Theorem~\ref{th:ccp-stationary}.The details are not repeated.
\end{proof}

\subsection{An illustrative example}

We use the following slight modification of the function (\ref{eq:deterministic 3-piece}) to illustrate the sign assumptions
in Theorems~\ref{th:fixed-point stationary} and \ref{th:ccp-stationary} on the simplified
problem: $\displaystyle{
\operatornamewithlimits{\mbox{\bf minimize}}_{x \in X}
} \, \Psi(x)$, where
\begin{equation} \label{eq:modified 3-piece}
\Psi(x) \, = \, \left\{ \begin{array}{ll}
\psi_1(x) & \mbox{if $a \leq f(x) < b$} \\ [5pt]
\psi_2(x) & \mbox{if $f(x) < a$} \\ [5pt]
\psi_3(x) & \mbox{if $f(x) \geq b$}.
\end{array} \right.
\end{equation}
Similar to the previous derivation, we can write
\[ 
\Psi(x) \, = \,
\psi_3(x) + \left( \, \psi_1(x) - \psi_3(x) \, \right) \, \onebld_{( \, 0,\infty \, )}( b - f(x) ) +
\left( \, \psi_2(x) - \psi_1(x) \, \right) \, \onebld_{( \, 0,\infty \, )}( a - f(x) ) .
\]
Given a vector $\bar{x} \in X$, there are five pseudo stationarity problems (\ref{eq:barx stationary problem}) at $\bar{x}$
for the problem $\displaystyle{
\operatornamewithlimits{\mbox{\bf minimize}}_{x \in X}
} \ \Psi(x)$ depending on the value of $f(\bar{x})$; these problems are

\gap

$\bullet $ [ $\displaystyle{
\operatornamewithlimits{\mbox{\bf minimize}}_{x \in X}
} \ \psi_1(x)$ ] if $f(\bar{x}) \in ( \, a,b \, )$;

\gap

$\bullet $ [ $\displaystyle{
\operatornamewithlimits{\mbox{\bf minimize}}_{x \in X}
} \ \psi_2(x)$ ] if $f(\bar{x}) \in ( \, -\infty, a \, )$;

\gap

$\bullet $ [ $\displaystyle{
\operatornamewithlimits{\mbox{\bf minimize}}_{x \in X}
} \ \psi_3(x)$ ] if $f(\bar{x}) \in ( \, b, \infty \, )$;

\gap

$\bullet $ [ $\displaystyle{
\operatornamewithlimits{\mbox{\bf minimize}}_{x \in X}
} \ \psi_1(x)$ {\bf subject to} $f(x) \geq a$ ] if $f(\bar{x}) = a$;

\gap

$\bullet $ [ $\displaystyle{
\operatornamewithlimits{\mbox{\bf minimize}}_{x \in X}
} \ \psi_3(x)$ {\bf subject to} $f(x) \geq b$ ] if $f(\bar{x}) = b$.

\gap

Clearly, not all stationary solutions of the above five problems are
local minimizers of $\Psi$ on $X$; nevertheless, they provide candidate minimizers.
The sign conditions in Theorem~\ref{th:fixed-point stationary} require that:

\gap

$\bullet $ $[ f(x) = a, x \in X ] \ \Rightarrow \ \psi_2(x) \geq \psi_1(x)$; and

\gap

$\bullet $ $[ f(x) = b, x \in X ] \ \Rightarrow \ \psi_1(x) \geq \psi_3(x)$.

\gap

Under these conditions, which basically stipulate that the function $\Psi$ can not rise up at a point of discontinuity,
the epi-hypographical approach will produce a pseudo B-stationary solution of
$\Psi$ on $X$.  In contrast, the pointwise sign stipulations at $x^*$ in condition (i) of Theorem~\ref{th:ccp-stationary} require that

\gap

$\bullet $ $f(x^*) = a \ \Rightarrow \ \psi_2(x) \geq \psi_1(x)$ in a neighborhood of $x^*$, and

\gap

$\bullet $ $f(x^*) = b \ \Rightarrow \ \psi_1(x) \geq \psi_3(x)$ in a neighborhood of $x^*$.

\gap

Condition (\ref{eq:strong convex gk}) becomes:

\gap

$\bullet $ there exists a neighborhood ${\cal N}$ of $\bar{x}$ such that
\[ \left. \begin{array}{l}
f^{\, \prime}(x^*;v) \, \geq \, 0 \\ [5pt]
v \, \in \, {\cal T}(X;x^*)
\end{array} \right\} \ \Rightarrow \ f^{\, \prime}(x;v) \, \geq \, 0 \epc \forall \, x \in {\cal N}.
\]
Along with the objective consistency condition of the iterates,
the above sign stipulations provide sufficient conditions for $x^*$ to be
a pseudo B-stationary solution of $\Psi$ on $X$.

\section{Concluding Remarks}

At the completion of the paper, the authors are grateful to receive a preprint \cite{Royset21} in which the author introduces
two concepts of consistent approximations in composite optimization; see Definition~2.2 therein.  The setting
of this reference
is the minimization of an objective that is the sum of an extended-valued constraint indicator function
and a composite function $h \circ F(x)$ where $h$ is a convex extended-valued function and $F$ is a vector
function whose components are locally Lipschitz continuous functions.  While in principle, the consistency concepts
defined therein
can be extended to any optimization problem without bother of its structure, the challenge is twofold: (a) how to define
the approximations, and (b) establishing the consistency of the derived approximations.  Although the problem
(\ref{eq:approximated POP}) is an approximation of (\ref{eq:original POP}), it is rather doubtful if the theory in
the reference could be applied to our context.  For one thing, a product $\psi \onebld_{( \, 0,\infty )}(f)$ is
very different from the composite family $h \circ F$ with $h$ and $F$ as stated.  As one can see, our analysis
makes extensive use of the product form and properties of the approximations of the Heavisde functions, in particular
requiring assumptions that are akin to such a structure.  Hopefully, our work will provide
a motivation to extend the theory of consistent approximations to broader composite classes of discontinuous functions.

\gap

{\bf Acknowledgements.}  The authors are grateful to the referees for their insightful comments that have helped to 
improve the presentation of the manuscript.

\gap

 {\bf Appendix: (C4) and weak sharp minima.} Consider the optimization problem
\begin{equation} \label{eq:opt for wsp}
\displaystyle{
\operatornamewithlimits{\mbox{\bf minimize}}_{x \in S}
} \ f(x),
\end{equation}
where $f$ is a continuous function bounded below on the closed set $S$.  Let
\[
f_{\min} \, \triangleq \displaystyle{
\operatornamewithlimits{\mbox{\bf minimum}}_{x \in S}
} \ f(x) \epc \mbox{and} \epc
{\cal F}_{\min} \, \triangleq \, \displaystyle{
\operatornamewithlimits{\mbox{\bf argmin}}_{x \in S}
} \ f(x).
\]
The problem (\ref{eq:opt for wsp}), or the pair $(f,S)$, is said to have weak sharp minima if
there exists a constant $\eta > 0$ such that
\[
f(x) - f_{\min} \, \geq \, \eta^{-1} \, \dist(x,{\cal F}_{\min}) \epc \forall \, x \, \in \, S.
\]
The definition of weak sharp minima was introduced in the Ph.D.\ thesis of Ferris \cite{Ferris88}.
Extensive discussion of this property and its role in optimization can be found in \cite[Section~6.5]{FacchineiPang03}.
The proposition below contains 3 statements.  Statement (a) asserts the existence of weak sharp minima for
(\ref{eq:opt for wsp}).  Statement (b) is known as Takahashi condition \cite{Takahashi91} in nonlinear analysis;
it implies in particular the existence of a global minimizer to the optimization
problem in question.  Statement (c) is the key to connecting condition (C4) to weak sharp minima.  Proof of
the implication (b) is by the renowned Ekeland's variational principle \cite{Ekeland74}; proof of the proposition
can be found in \cite[Section~6.5]{FacchineiPang03}; see also \cite[Section~8.5.4]{CuiPang2021}.

\begin{proposition} \label{pr:takahashi} \rm
Let $S$ be a closed set in $\mathbb{R}^n$ and $f$ be a continuous
real-valued function defined and bounded below on $S$.  Let $f_{\rm inf}$ denote
the infimum value of $f$ on $S$.
The two statements (a) and (b) are equivalent and either
statement is implied by statement (c):

\gap

\noindent (a) The pair $(f,S)$ has weak sharp minima; i.e., ${\cal F}_{\min} \neq \emptyset$
and for some constant $c > 0,$
\begin{equation} \label{eq:weak sharp general fnc}
f(x) - f_{\min} \, \geq \, c \, \dist(x,{\cal F}_{\min}), \epc \forall \, x \, \in \, S.
\end{equation}
(b) There exists a constant $c > 0$ such that for each $x \in S$ with $f(x) > f_{\rm inf}$, there
exists a vector $y \in S$ distinct from $x$ such that
\begin{equation} \label{eq:takahashi inequality}
f(y) + c \, \| \, x - y \, \| \, \leq \, f(x).
\end{equation}
(c) Suppose that $f$ is B-differentiable on $S$ and
there exists a scalar $\delta > 0$ such that
for every $x \in S$ with $f(x) > f_{\rm inf}$, a vector
$d \in {\cal T}(S;x)$ with unit (Euclidean) length exists
satisfying $f^{\, \prime}(x;d) \leq -\delta$.

\gap

\noindent Moreover, if $f$ is a convex function and $S$ is a convex set, then (b) implies (c);
so all three statements (a), (b), and (c) are equivalent. \hfill $\Box$
\end{proposition}

To apply Proposition~\ref{pr:takahashi} to condition (C4), we rewrite the various index sets in $\wh{S}_{\rm ps}(x)$
more generally.  Specifically, for any two triplets of index sets $\boldsymbol{\cal K} \triangleq \{ {\cal K}_>, {\cal K}_=, {\cal K}_< \}$
and $\boldsymbol{\cal L} \triangleq \{ {\cal L}_>, {\cal L}_=, {\cal L}_< \}$ with the former partitioning $\{ 1, \cdots, K \}$ and
the latter partitioning $\{ 1, \cdots, L \}$, define the closed set
\[
{\cal S}(\boldsymbol{\cal K},\boldsymbol{\cal L}) \, \triangleq \, \left\{ \, x \, \in \, X \, \left| \, \begin{array}{ll}
g_k(x) \, \leq \, 0 & \forall \, k \, \in \, {\cal K}_< \, \cup \, {\cal K}_= \\ [0.1in]
g_k(x) \, \geq \, 0 & \forall \, k \, \in \, {\cal K}_> \\ [0.1in]
h_{\ell}(x) \, \leq \, 0 & \forall \, \ell \, \in \, {\cal L}_< \, \cup \, {\cal L}_= \\ [0.1in]
h_{\ell}(x) \, \geq \, 0 & \forall \, k \, \in \, {\cal L}_>
\end{array} \right. \, \right\}.
\]
and consider the optimization problem:
\[ \displaystyle{
\operatornamewithlimits{\mbox{\bf minimize}}_{x \in {\cal S}(\boldsymbol{\cal K},\boldsymbol{\cal L})}
} \ f_{\boldsymbol{\cal K},\boldsymbol{\cal L}}(x) \, \triangleq \, \max\left( \, \displaystyle{
\sum_{\ell \in {\cal L}_>}
} \, \phi_{\ell}(x) - b, \, 0 \, \right).
\]
Assume that each function $\phi_{\ell}$ is bounded below on $X$.  Then a very loose sufficient condition for (C4) to hold
is that for all pairs $(\boldsymbol{\cal K},\boldsymbol{\cal L})$, the minimum value of $f_{\boldsymbol{\cal K},\boldsymbol{\cal L}}$
on ${\cal S}(\boldsymbol{\cal K},\boldsymbol{\cal L})$ is zero and the pair
$\left( \, f_{\boldsymbol{\cal K},\boldsymbol{\cal L}}, \, {\cal S}(\boldsymbol{\cal K},\boldsymbol{\cal L}) \, \right)$
satisfies condition (c) of Proposition~\ref{pr:takahashi} with $\delta = 1$.  We believe that it may be possible to tighten this sufficient
condition and derive a result similar to the proposition for condition (C4); such details are beyond the scope of this work and best left
for a separate investigation.

\end{document}